\documentclass[a4paper]{article}

\usepackage[a4paper,
            bindingoffset=0.2in,
            left=1in,
            right=1in,
            top=1in,
            bottom=1in,
            footskip=.25in]{geometry}

\usepackage{graphicx}
\usepackage{float}
\usepackage{subfigure}
\usepackage{amsfonts, amsmath, amsthm, amssymb}
\usepackage{authblk}
\usepackage{tikz}
\usepackage{pgfplots}
\pgfplotsset{compat=newest}
\usepackage[maxbibnames=9]{biblatex} 
\usepackage{hyperref}
\usepackage{comment}
\usepackage{mathtools}

\usepackage{algorithm}
\usepackage{algpseudocode}

\newtheorem{theorem}{Theorem}[section]
\newtheorem{lemma}[theorem]{Lemma}
\newtheorem{definition}[theorem]{Definition}
\newtheorem{remark}[theorem]{Remark}
\newtheorem{proposition}[theorem]{Proposition}
\newtheorem*{nat}{Notation and terminology}

\def\R{{\mathbb R}}

\def\N{{\mathbb N}}

\def\P{{\mathbb P}}

\newcommand{\de}{\,\mathrm{d}}
\newcommand{\dist}{\,\mathrm{dist}}
\newcommand{\dx}{\,\mathrm{d} x_1 \ldots \mathrm{d} x_n}

\newcommand{\px}{p(x_1, \ldots, x_n)}

	{\left\lbrace \begin{array}{@{} l @{} }}%
	{ \end{array}\right.}

\ifpdf
  \DeclareGraphicsExtensions{.eps,.pdf,.png,.jpg}
\else
  \DeclareGraphicsExtensions{.eps}
\fi

\usepackage{xcolor}
\definecolor{turchese}{RGB}{35, 174, 163}

\title{Interpolation by integrals on balls} 

\author[a,b]{Ludovico Bruni Bruno}
\author[a]{Giacomo Elefante}
\affil[a]{Department of Mathematics \lq\lq Tullio Levi-Civita\rq\rq, University of Padova, Via Trieste 63, Padova, Italy}
\affil[b]{Istituto Nazionale di Alta Matematica \lq\lq Francesco Severi\rq\rq, Piazzale Aldo Moro 5, Roma, Italy}
\affil[ ]{\texttt{ludovico.brunibruno@unipd.it}, \texttt{giacomo.elefante@unipd.it}}
\date{}

\begin{document}

\maketitle

\begin{abstract}
In this work we blend interpolation theory with numerical integration, constructing an interpolator based on integrals over $n$-dimensional balls. We show that, under hypotheses on the radius of the $n$-balls, the problem can be treated as an interpolation problem both on a collection of $(n-1)$-spheres $ S^{n-1} $ and multivariate point sets, for which a wide literature is available. With the aim of exact quadrature and cubature formulae, we offer a neat strategy for the exact computation of the Vandermonde matrix of the problem and propose a meaningful Lebesgue constant. Problematic situations are evidenced and a charming aspect is enlightened: the majority of the theoretical results only deal with the centre of the domains of integration and are not really sensitive to their radius. We flank our theoretical results by a large amount of comprehensive numerical examples.
\end{abstract}

\section{Introduction}

The principal aim of this work is answering the following question: let $ p(x_1, \ldots, x_n) $ be a $n$-variate polynomial and let $ B^n ({\bf c}, R) $ be the $n$-dimensional sphere of radius $ R>0 $ centered at $ {\bf c} \in \R^n $
$$ B^n ({\bf c}, R) \coloneqq \left\{ {\bf x} \in \R^n \ \mid \ \Vert {\bf x} - {\bf c} \Vert_2 \leq R \right\} .$$
How can we place a collection of balls $ B^n_i ({\bf c}_i, r_i) \subseteq B^n ({\bf c}, R) $ such that
$$ \int_{B^n_i ({\bf c}_i, r_i)} p(x_1, \ldots, x_n) \de x_1 \ldots \de x_n = 0 \hspace{.2cm} \forall i \quad \text{implies that} \quad p(x_1, \ldots, x_n) = 0 \ ?$$
To provide an answer to this problem we shall transform it into an interpolation problem. 
The solution of the problem may involve lattices on the $n$-simplex
$$ T \coloneqq \left\{ {\bf x} \in \R^n \ \mid \ {\bf x} = \sum_{i=0}^n \lambda_i {\bf v}_i, \ 0 \leq \lambda_i \leq 1, \ \det[{\bf v}_1 - {\bf v}_0 | \ldots | {\bf v}_n - {\bf v}_0] \ne 0  \right\} $$
as well as more general scattered points in the $n$-ball, but it may also concern orbits under the action of the group $ SO(n) $.
All these settings make it possible to appeal a consolidated literature. For the case of simplices we refer to \cite{BruniThesis}, where the problem is already treated in terms of integrals, as well as on the more classical literature of multivariate interpolation on simplices \cite{BSV12,Hesthaven98,Warburton06} and polygons \cite{GSV11}. For the case of discs and balls (possibly $n$-dimensional) we refer to \cite{MS19,Phung21,Phung17,Xu03,Xu04}. These two approaches coincide when $ n = 1 $; this has been treated from different perspectives \cite{Bojanov06,BE23}.

The resulting interpolation operator thus associates intregrable functions with multivariate polynomials that have prescribed integrals. The norm of the interpolation operator can be characterised in a charming geometrical way that makes it possible to perform computations and numerical tests.

In stark contrast with what has been observed in the simplicial case (see, e.g., \cite{ABR23,BruniThesis}), where both the positioning and the size of the supports of integration play a relevant role, in this work we are able to separate these two aspects: the centres of the balls (hence their placement) carry the majority of the interpolation features, both on the theoretical side and the numerical one. In fact, we shall prove that it is only the centres of the balls that determine the unisolvence of the desired set. The role of the radius become significant only when one wants to characterise the norm of the interpolation operator as an extension to this framework of the generalised Lebesgue constant given in \cite{ARR20}. Interestingly, numerical tests show that the quality of the interpolation is again significantly more affected by the centres than by the radii of the supporting balls.

Although we do not pursue this path here yet, we claim that the techniques we present can be extended to any domain which can be obtained via affine mappings from $ B^n ({\bf c}, R) $, e.g. ellipsoids. In addition, results on unisolvence turn out to be effective also on very general domains, independently of their shape. We reckon that this might provide a first step into blending techniques of weights \cite{Rapetti07,BruniThesis} with kernel methods \cite{BuhmannIntegration}, such as radial basis functions (RBF) methods \cite{buhmann_2003}.

\textbf{Outline of the paper}. The paper is organised as follows. In Section \ref{sect:theory} we provide theoretical results that answer the main question of the paper. Section \ref{sect:interpolationtheory} deals with the features of the interpolator that naturally arises from the above question. A concept of quality is introduced and characterised in terms of a meaningful Lebesgue constant. From Section \ref{sect:r2} we find it convenient to fix the dimension of the ambient space to two, in order to present a more visualisable description. In view of this, we are able to offer in Section \ref{sect:fails} an overview of relevant fails. The numerical side is condensed in Section \ref{sect:numeric}, where a distinct subsection is dedicated to each of the following topics: conditioning and choice of the basis, interpolation and stability (in terms of the aforementioned Lebesgue constant). Conclusions are gathered at the end of the paper, in Section \ref{sect:conclusions}.

\section{Two approaches to unisolvence} \label{sect:theory}

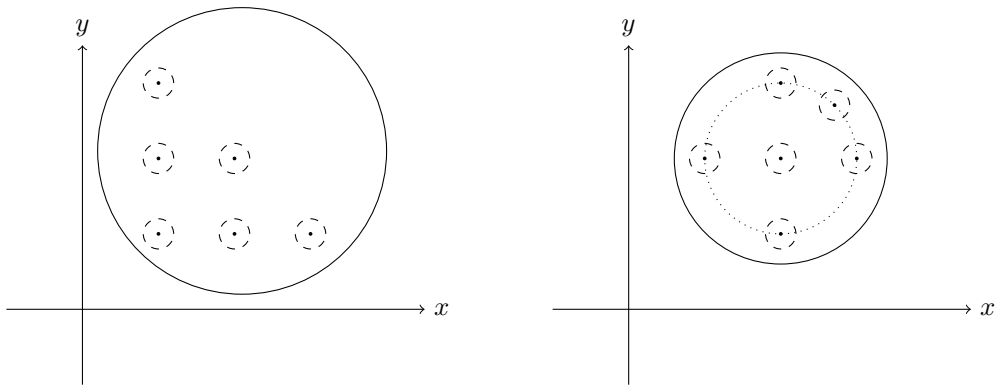
\begin{figure}[H]
\centering
\begin{tikzpicture}
  \draw[->] (-1, 0) -- (4.5, 0) node[right] {$x$};
  \draw[->] (0, -1) -- (0, 3.5) node[above] {$y$};
  \filldraw[black] (1,1) circle (0.5pt);
  \filldraw[black] (2,1) circle (0.5pt);
  \filldraw[black] (3,1) circle (0.5pt);
  \filldraw[black] (1,2) circle (0.5pt);
  \filldraw[black] (2,2) circle (0.5pt);
  \filldraw[black] (1,3) circle (0.5pt);
  
\draw[thin,dashed] (1,1) circle (0.2cm);
\draw[thin,dashed] (2,1) circle (0.2cm);
\draw[thin,dashed] (3,1) circle (0.2cm);
\draw[thin,dashed] (1,2) circle (0.2cm);
\draw[thin,dashed] (2,2) circle (0.2cm);
\draw[thin,dashed] (1,3) circle (0.2cm);

\draw[thin] (2.1,2.1) circle (1.9cm);
\end{tikzpicture} \hspace{1cm}
\begin{tikzpicture}
  \draw[->] (-1, 0) -- (4.5, 0) node[right] {$x$};
  \draw[->] (0, -1) -- (0, 3.5) node[above] {$y$};
  \filldraw[black] (3,2) circle (0.5pt);
  \filldraw[black] (2,3) circle (0.5pt);
  \filldraw[black] (1,2) circle (0.5pt);
  \filldraw[black] (2,1) circle (0.5pt);
  \filldraw[black] (2,2) circle (0.5pt);
  \filldraw[black] (2.7071,2.7071) circle (0.5pt);
  
\draw[thin,dashed] (3,2) circle (0.2cm);
\draw[thin,dashed] (2,3) circle (0.2cm);
\draw[thin,dashed] (2,1) circle (0.2cm);
\draw[thin,dashed] (1,2) circle (0.2cm);
\draw[thin,dashed] (2,2) circle (0.2cm);
\draw[thin,dashed] (2.7071,2.7071) circle (0.2cm);

\draw[thin] (2,2) circle (1.4cm);
\draw[thin,dotted] (2,2) circle (1cm);
\end{tikzpicture}
\label{fig:variousunisolvence}
\caption{Integration on discs is unisolvent for $ \mathbb{P}_2 (\R^2)$ both on the left (it is a principal lattice) and on the right (5 points lie on the same orbit and the last one is outside).} 
\end{figure}

Theoretical results of this section follow two directions: the first enlightens the role of lattices and points sets, as in interpolation over simplices or cubes, and is shown in Figure \ref{fig:variousunisolvence}, left; the second exploits the regularity of the shape taken into account, whence the role of orbits, and is shown in Figure \ref{fig:variousunisolvence}, right. Both techniques will turn out to be blind to the radius $ R $ of $ B^n ({\bf c}, R) $.

\begin{nat}
    In the following we shall frequently adopt the (not restrictive) simplification in which $ {\bf c} = {\bf 0 }$ and $ R = 1 $, which justifies the notation $ B^n \coloneqq B({\bf 0}, 1) $ whence $ \partial B^n = S^{n-1} $. Formally, this is equivalent to produce the change of variable
$$ (x_1, \ldots, x_n) \mapsto \frac{1}{R} (x_1 - c_1 , \ldots, x_n - c_n) ,$$
which is affine and hence does not affect polynomial objects.
\end{nat}

\subsection{Unisolvence by points} \label{sect:unisolvencepoints}

To begin with, we consider the case that can be reduced to unisolvence on point sets and lattices. The following lemma is a generalisation of \cite[Lemma $ 3.1 $]{BBZ22} and \cite[Lemma $3.12$]{ChR16} to this framework.

\begin{lemma} \label{lem:alltranslated}
    Let $ r > 0 $ and let $ B_{{\bf \xi}}^n \coloneqq B^n ({\bf \xi}, r) $. If
    \begin{equation} \label{eq:genericxi}
        \int_{B_{{\bf \xi}}^n} \px \de x_1 \ldots \de x_n = 0
    \end{equation}
    for each $ \xi \in \R^n $, then $ \px = 0 $.
\end{lemma}

\begin{proof}
    Let us denote by $ \mathcal{Z}[\px] $ the zero locus of $ \px $, which is an algebraic variety. Since the measure of $\mathcal{Z}[\px]$ is equal to zero, for each $ r > 0 $ there exists $ \xi \in \mathbb{R}^n $ such that $ B_{\xi}^n \cap \mathcal{Z}[\px] = \emptyset $. As a consequence, $ \mathrm{sgn} \left( \px \right) $ is constant on $ B_{\xi}^n $ and its interior $ \mathring{B}_{\xi}^n $, which is open and non empty. Thus Eq. \eqref{eq:genericxi} implies that $ p({\bf x}) = 0 $ for each $ {\bf x} \in \mathring{B}_{\xi}^n $ and hence $ p({\bf x}) = 0 $  for each $ {\bf x} \in \mathbb{R}^n $.
\end{proof}

Lemma \ref{lem:alltranslated} is the key tool of all the subsequent techniques, as we will be able to place ourselves in the comfortable situation where Eq. \eqref{eq:genericxi} is satisfied. The following theorem evidences the relationship between unisolvent sets in nodal interpolation and the centres of the spheres $ B_{\xi}^n $.

\begin{theorem} \label{thm:generalunisolvence}
    Let $ \Xi \coloneqq \{ \xi_i \}_{i=1}^N $ be a unisolvent set for $ \mathbb{P}_d (\mathbb{R}^n) $ and $ r>0 $. If
    $$ \int_{B_{\xi_i}^n} \px \de x_1 \ldots \de x_n = 0 $$
    for any $ B_{{\bf \xi}_i}^n \coloneqq B^n ({{\bf \xi}_i},r) $, then $ \px = 0 $.
\end{theorem}
\begin{proof}
    Write each $ \{ B_{\xi_i}^n \}_{i=1}^N $ as the image of $ B^n$ via the affinity 
    $$ {\bf x} \mapsto r {\bf x} + {\bf \xi } .
    $$
    By this change of variable, it is immediate to recognise that the map
    $$ {\bf \xi} \mapsto q ({\bf \xi}) \coloneqq \int_{B_{\xi}^n} \px \de x_1 \ldots \de x_n $$
    is a polynomial of degree $ d = \deg \px $ in the variable $ \xi $. By hypothesis, $ q(\xi_i) = 0 $ for any $ \xi_i \in \Xi $. Since this set is unisolvent for $ \mathbb{P}_d $, it follows that $ q(\xi) $ is in fact the zero polynomial, hence $ q(\xi) = \int_{B_{\xi}^n} \px = 0 $ for each $ \xi \in \mathbb{R}^n $. Lemma \ref{lem:alltranslated} ensures that $ \px = 0 $.
\end{proof}

This construction makes it possible to appeal to a large choice of unisolvent sets that are known on simplicial or tensor product domains. To account some of them, we refer to \cite{Bos91,BSV12,GSV11,Hesthaven98}.

\begin{remark}
    Constraints on $ r $ shall be taken into account only when interested in sets that do not exceed the boundary $ \partial B^n ({\bf c}, R) $ of $ B^n ({\bf c}, R) $, for instance when dealing with neighbouring elements. If we set 
    $ \rho \coloneqq \dist(\Xi, \partial B^n ({\bf c}, R) ) $, we immediately retrieve the condition $ 0 < r \leq \rho $.
\end{remark}

\subsection{Unisolvence by orbits} \label{sect:unisolvenceorbits}

The structure of $ B^n $ allows for a neat identification of unisolvent sets $ \Xi $, as pointed out in \cite{Bos91}. We closely follow such a construction to extend it to our framework. From now on the simplification on the reference domain will be systematic.

The algebraic variety
$$ S_r^{n-1} \coloneqq \left\{ (x_1, \ldots, x_n) \in \R^n \ \mid \ \sum_{i=1}^n x_i^2 = r^2 \right\} $$
describes the orbit of a point at distance $ r $ from $ {\bf 0} $ under the action of $ SO(n) $. It is supported in $ B^n $, provided that $ 0 < r \leq 1 $. Being orbits, they do not intersect. Moreover, since the polynomial 
$$ P(x_1, \ldots, x_n ) \coloneqq \sum_{i=1}^n x_i^2 - r^2 $$
that describes $ S_r^{n-1} $ is irreducible, the ideal associated with $ P(x_1, \ldots, x_n ) $ is principal. Thus \cite[Lemma $2.1$]{Bos91} prescribes the dimension of the space of polynomials on $ S_r^{n-1} $ to be
\begin{align} \label{eq:dimrestr}
    \dim \left( \mathbb{P}_d \left( S_r^{n-1} \right) \right) & = \dim \left( \mathbb{P}_d \left( \R^n \right) \right) - \dim \left( \mathbb{P}_{d-2} \left( \R^n \right) \right) \nonumber \\
    & = \binom{d+n}{n} - \binom{d+n-2}{n},
\end{align}
where we have used the fact that $ \deg P(x_1, \ldots, x_n) = 2 $. 

\begin{definition}
We say that an orbit $ S_r^{n-1} $ is \emph{saturated} for $ \P_d (S_r^{n-1}) $ if it contains $ \dim \P_d (S_r^{n-1}) $ points.\end{definition}
%

The concept of saturation is relevant in the following description. In fact, saturated orbits prescribe the largest number of supports of integration that can be centered on $ S_r^{n-1} $. Of course saturation, as well as the dimension count \eqref{eq:dimrestr}, are actually independent on the radius $r$; its inclusion in the notation is useful for a better identification of the orbit.

\begin{lemma} \label{lem:factorisationSr} 
Let $ p(x_1, \ldots, x_n) \in \P_d (\R^n) $ and let $ S_r^{n-1} $ be saturated for $ \P_d (S_r^{n-1}) $. Suppose $ \Xi $ is unisolvent for $ \mathbb{P}_d \left( S_r^{n-1} \right) $. Then, if
$$ \int_{B_{\xi}^r} p(x_1, \ldots, x_n) = 0 $$
for each $ \xi \in \Xi $, $ P(x_1, \ldots, x_n) $ divides $ \px $.
\end{lemma}

\begin{proof}
    Since $ \Xi $ is unisolvent for $ S_r^{n-1} $, which is a saturated orbit, by Theorem \ref{thm:generalunisolvence} we get that
    $$ \int_{B_{\xi}^r} p(x_1, \ldots, x_n) = 0 $$
    for each $ \xi \in \Xi $ implies that $ p(x_1, \ldots, x_n) \bigl|_{S_r^{n-1}} = 0 $, hence $ p(x_1, \ldots, x_n) = P(x_1, \ldots, x_n) q(x_1, \ldots, x_n) $ for some $ q(x_1, \ldots, x_n) \in \P_{d-2}(\R^n) $.
\end{proof}

With Lemma \ref{lem:factorisationSr} in mind, we produce the following Algorithm \ref{alg:orbits}, whose well posedness will be a consequence of Theorem \ref{thm:unisolvenceorbits}.

\begin{algorithm}
\caption{Construction of the unisolvent set by orbits}\label{alg:orbits}
\begin{algorithmic}
\Require The total degree $ d $
\Ensure A collection of unisolvent balls
\State $ {\bf r} = [ ] $ \Comment{Vector of the radii}
\While{$d \geq 1 $}
\State choose $ r \not\in {\bf r} $ and add $ r $ to $ {\bf r }$ \Comment{Guarantees disjointess of orbits}
\State fix $ r_j $ \Comment{Radius of the balls on this orbit}
\State compute $ d_j \coloneqq \dim(\P_d (S^{n-1}_r)) $ by \eqref{eq:dimrestr} \Comment{Number of balls to place on the $j$-th orbit}
\State choose $c_j$ which are $ d_j $ points on the orbit $ S_{r_j} $
\State construct $ d_j $ balls of radius $ r_j $ centred on $ c_j $
\State $ d = d-2 $ \Comment{Every saturated orbit reduces by $ 2 $ the total degree}
\EndWhile
\If{$ d = 0$} \Comment{Only the degenerate orbit is left}
\State construct a ball of any radius centred at $ {\bf 0}$
\EndIf 
\end{algorithmic}
\end{algorithm}
An iterative application of Lemma \ref{lem:factorisationSr} on distinct orbits is the key tool of the above algorithm and the core of the proof of the next result.
%

\begin{theorem} \label{thm:unisolvenceorbits} Let $ \px \in \P_d (\R^n) $. Let $ \Xi $ be a collection of points such that $ \Xi_j \coloneqq S_{r_j}^{n-1} \cap \Xi $ is unisolvent for $ \P_{d_j} ( S_{r_j}^{n-1} ) $, with $ d_j = 0, 2, \ldots, d $ if $ d $ is even or $ d_j = 1, 3, \ldots, d $ if $ d $ is odd. If
    $$ \int_{B^n ({\xi_j}_i, r_j)} \px \dx = 0 $$
    for each $ j_i $, then $ \px = 0 $.
\end{theorem}

\begin{proof}
    Call $ d^* $ the number of saturated orbits. 
    For $ j = 1, \ldots, d^* $, order orbits $ \Xi_j $ so that $ j > k $ if $ \# \left( \Xi_j \right) > \# \left( \Xi_k \right) $, the symbol $ \# $ denoting the cardinality of centres lying on that orbit. Since
    $$ \int_{B^n ({\xi_j}_i, r_j)} \px \dx = 0 $$
    for each ball, 
    in particular,
    $$ \int_{B^n ({\xi_1}_i, r_1)} \px \dx = 0 $$
    for the first orbit, i.e., for $ j = 1 $ and for each $ i $. Hence, by Lemma \ref{lem:factorisationSr}, there exists a degree $ 2 $ polynomial
    $$ P_1 (x_1, \ldots, x_n) \coloneqq \sum_{i=1}^n x_i^2 - r_1^2 $$
    such that $ \px = P_1 (x_1, \dots, x_n) q_1 (x_1, \ldots, x_n) $, with $ \deg q_1 = d-2 $. 
    Iterating this up to $ j = d^* -1 $, we get
    $$ \px = \left( \prod_{j=1}^{d^* - 1} P_j (x_1, \ldots, x_n) \right) q_{d^* -1} (x_1, \ldots, x_n) ,$$
    with $ q_{d^* -1} $ a polynomial of degree $ \bar{d} \leq 1 $. When restricting to the last orbit 
    $ S^{n-1}_{r_{d^*}} $, we have
    $$ \px \bigl|_{S^{n-1}_{r_{d^*}}} = c q_{d^*-1} (x_1, \ldots, x_n) ,$$
    with $ c = \prod_{j=1}^{d^*-1} (r_j - r_{d^*})^2 $ being a constant. We show, separating the cases $ \bar{d} = 0, 1 $, that     
    $$ \int_{B^n ({\xi_{d^*}}_i, r_{d^*})} q_{d^*-1}(x_1, \ldots, x_n) \dx = 0 $$
    for all $ {\xi_{d^*}}_i \in \Xi_{d^*} $ implies that $ q_{d^*-1}(x_1, \ldots, x_n) = 0 $. This will complete the proof.
    
    Let us first consider the case when $ \bar{d} = 0 $. The polynomial $ q_{d^*-1} $ is therefore a constant, say $ c' $, 
    with vanishing integral hence $ c' = 0 $ and $ \px = 0 $. 
    
    If we assume $ \bar{d} = 1 $ the claim is a consequence of Lemma \ref{lem:factorisationSr}. In fact, the algebraic variety described by $ P_{d^*}(x_1, \ldots, x_n) $ has degree $ 2 $ and it may not divide $ q_{d^*-1}(x_1, \ldots, x_n) $. Hence $ q_{d^*-1}(x_1, \ldots, x_n) = 0 $ and once again $ \px = 0 $.
\end{proof}

Note that the hypotheses on the radii in Theorem \ref{thm:unisolvenceorbits} are weaker than those of Theorem \ref{thm:generalunisolvence}: indeed, a different fixed radius $ r_j $ for each orbit is allowed. Moreover, if all the radii are equal, this offers a different proof of Theorem \ref{thm:generalunisolvence}, without needing a global unisolvence result. On the contrary, when orbits are not saturated, a global unisolvence result will be needed to invoke Theorem \ref{thm:generalunisolvence}; the necessity of this hypothesis is studied in Section \ref{sect:fail2}.

\section{Interpolation and Lebesgue constant} \label{sect:interpolationtheory}

As long as unisolvence is achieved, we may easily define an interpolator by asking a polynomial and a function to share the same degrees of freedom. In particular, $\px$ is the interpolated of $ f(x_1, \ldots, x_n) $ if
\begin{equation} \label{eq:interpolatingconditions}
    \int_{B_i} \px \dx  = \int_{B_i} f(x_1, \ldots, x_n) \dx 
\end{equation}
for each $ i $. Unisolvence guarantees the well posedness of the interpolator associated with \eqref{eq:interpolatingconditions}. A useful representation is then that obtained by exploiting Lagrange basis functions $ \ell_{B_i} $ such that
\begin{equation} \label{eq:dualityLB}
    \int_{B_j} \ell_{B_i} = \delta_{i,j} .
\end{equation}
Explicit expressions and closed formulae for \eqref{eq:dualityLB} are not always at hand. We, therefore, proceed as follows. 

Let us set $ N \coloneqq \dim \mathbb{P}_d (\R^n) $. We recall that the unisolvence of the collection $ \{ B_i \}_{i=1}^N $ does not depend on the basis of $ \mathbb{P}_d (\R^n) $ chosen. Once a convenient basis $ \{ p_i \}_{i=1}^N $ is fixed, we may construct the corresponding Vandermonde matrix
\begin{equation} \label{eq:VdM}
V_{i,j} \coloneqq \int_{B_j} \px \dx . 
\end{equation}
The coefficients of the Lagrange basis $ \{ \ell_{B_i} \}_{i=1}^N $ respect to the basis $ \{ p_i \}_{i=1}^N $ are constructed through the inverse of this matrix. In fact, it is immediate to check that 
\begin{equation} \label{eq:LagrBases}
\ell_{B_i} \coloneqq \sum_{j=1}^N V^{-1}_{i,j} p_i ,
\end{equation}
satisfies \eqref{eq:dualityLB} for each $i=1,\dots,N$.

\begin{remark}
    Since supports $ B_j $ are $n$-balls, elements of \eqref{eq:VdM} may be exactly computed by means of exact quadrature (cubature, and so on) formulae; see for the disc, e.g., \cite{CK00,TFR22} and references therein.
\end{remark}

With the aim of the Lagrange basis, it is immediate to produce the interpolator 
\begin{align} \label{eq:definterp}
    \Pi: \quad C^0 (\R^n) & \longrightarrow \P_d (\R^n) \nonumber \\
    f & \longmapsto \Pi f \coloneqq \sum_{i=1}^N \left( \int_{B_i} f \right) \ell_{B_i}.
\end{align}
Note that the regularity of the interpolated function can be weakened.

\begin{proposition} \label{prop:projector}
    Equation \eqref{eq:definterp} defines a projector onto $ \P_d (\R^n) $.
\end{proposition}

\begin{proof}
    By direct computation, applying \eqref{eq:dualityLB}:
    \begin{align*}
        \Pi \left(\Pi f\right) & = \sum_{j=1}^N \int_{B_j} \left( \sum_{i=1}^N \left( \int_{B_i} f \right) \ell_{B_i} \right) \ell_{B_j} = \sum_{j=1}^N \left( \sum_{i=1}^N \left( \int_{B_i} f \right)  \int_{B_j} \ell_{B_i} \right) \ell_{B_j} \\
        & = \sum_{j=1}^N \left( \int_{B_j} f \right) \ell_{B_j} = \Pi f,
    \end{align*}
    whence the claim.
\end{proof}

\subsection{The Lebesgue constant}
In classical interpolation problems, the Lebesgue constant is important because it expresses a way to compute how close the interpolant is to the best approximant and it may also provide the conditioning of the interpolant. For these reasons, many studies are devoted to find optimal points for the growth of the Lebesgue constant or study bounds for the growth of the interpolant at good set of points (e.g., \cite{BE21,DEM21,Ibrahimoglu16,MS19}).

Whenever a norm is placed on the domain and the codomain of $ \Pi $, the operator norm $ \Vert \Pi \Vert_{\mathrm{op}} $ is given. Following \cite{ARR20} and \cite{BruniThesis}, we find it convenient to consider the so-called $0$-norm studied in \cite{HarrisonLeb}
\begin{equation} \label{eq:zeronorm}
\Vert f \Vert_0 \coloneqq \sup_{c \in \mathcal{C}^n (\Omega) } \frac{1}{|c|_0} \left\vert \int_{c} f(x_1, \ldots, x_n) \dx \right\vert .
\end{equation}
The supremum in \eqref{eq:zeronorm} is sought over the collection of $n$-chains $ \mathcal{C}^n $ supported in some domain $ \Omega $ and appropriately measured. This perfectly matches the 
classical setting of simplicial \emph{weights} \cite{Rapetti07}, which is that considered in 
\cite{ARR20}. 
As we leave the simplicial setting, we modify the norm \eqref{eq:zeronorm} in two directions:
\begin{itemize}
    \item[i)] {\it shape of the supports}: we consider spheres in place of simplices, to make them coherent with subsequent computations; 
    \item[ii)] {\it connectedness of supports}: we remove the possibility of considering chains. If, on the one hand, the construction in \cite{HarrisonLeb} and the definition in \cite{ARR20} allow for very developed supports, in practice the resulting quantity turns out to be very hard to be handled (see, e.g., the computations in \cite{BruniThesis}). 
    Some numerical results in the simplicial setting suggest that this choice is fine enough to capture the numerical outline of the problem \cite{ABR20}.
\end{itemize}
We thus design the norm
\begin{equation} \label{eq:ballnorm}
\Vert f \Vert_B \coloneqq \sup_{B \in \mathcal{B}^n (\Omega) } \frac{1}{|B|} \left\vert \int_{B} f(x_1, \ldots, x_n) \dx \right\vert .
\end{equation}
Here the space $ \mathcal{B}^n $ is that of all possible $n$-balls contained in $ \Omega $ and
\begin{equation} \label{eq:ballvolume}
    | B | = \frac{\pi^n}{\Gamma(\frac{n}{2}+1)} R^n
\end{equation}
is the volume of the $n$-ball of radius $R$. 

\begin{lemma}
    The quantity \eqref{eq:ballnorm} defines a seminorm on locally integrable functions on $ \Omega $. It defines a norm on continuous functions on $ \Omega $.
\end{lemma}

\begin{proof}
Axioms of seminorms are immediately verified as soon as $ \int_B f $ is meaningful.

We now add the hypothesis of continuity of $ f $ and show that, in this case, \eqref{eq:ballnorm} is also a norm. Let $ f \in C^0 (\Omega) $ and suppose, by contradiction, that there exists an $ f \ne 0 $ such that $ \Vert f \Vert_B = 0 $. Consider $ {\bf \xi} \in \Omega $ such that $ f({\bf \xi}) \ne 0$, say for example $ f({\bf \xi}) > 0 $. Hence there exists an open, non empty set $ U_\xi \subset \Omega $ such that $ \xi \in U_{\bf \xi} $ where $ f({\bf x}) > 0 $ for each $ {\bf x} \in U_{\bf \xi} $. Thus, for any ball $ B_{\bf \xi} \subset U_{\bf \xi} $, we get
$$ \left\vert \int_{B_{\bf \xi}} f(x_1, \ldots, x_n) \dx \right\vert \geq \min \left\vert f \bigl|_{B_{\bf \xi}} ({\bf x}) \right\vert \left\vert B_{\bf \xi} \right\vert > 0 ,$$
since $ B_{\bf \xi} $ is compact. This is in contradiction with the hypothesis and, therefore, the result is established.
\end{proof}

Endowing $ C^0 (\R^n) $ and $ \P_d (\R^n) $ with the norm \eqref{eq:ballnorm}, we get that the Lebesgue constant of the problem is
\begin{equation} \label{eq:defgeneralLeb}
    \Lambda_d \coloneqq \Vert \Pi \Vert_{\mathrm{op}} = \sup_{f \in C^0 (\R^n)}\frac{\Vert \Pi f \Vert_B}{\Vert f \Vert_B} .
\end{equation}

Keeping \eqref{eq:defgeneralLeb} under control is essential for having a reliable interpolation. As a consequence, the possibility of computing and estimating such a quantity is crucial. We provide a characterisation of \eqref{eq:defgeneralLeb} under appropriate assumptions.

\begin{theorem}
    Let $ \{ B_i \}_{i=1}^N $ be a collection of unisolvent $n$-balls for $ \P_d (\R^n) $. If $ B_i \cap B_j = \emptyset $ for $ i \ne j $, then
    \begin{equation} \label{eq:newLeb}
        \Lambda_d = \sup_{B \in \mathcal{B}^n (\Omega)} \frac{1}{|B|} \sum_{i=1}^N | B_i | \left\vert \int_B \ell_{B_i} \right\vert .
    \end{equation}
\end{theorem}

\begin{proof}
    We construct a sequence of functions $ \{ f_n \}_{n \in \N} \in C^0 (\R^n) $ such that 
    $$ \Vert \Pi \Vert_{\mathrm{op}} \geq \Lambda_d - \frac{1}{r' n} ,$$
    with $ r'$ being a constant independent of $n$. Define $ B_i^{(n)} \coloneqq B(\xi_i, r_i - \frac{1}{n}) $ and, for each $ B_i^{(n)} $, the function
    $$ f_i^{(n)} ({\bf x}) =
    \begin{cases}
    1 \qquad\qquad\qquad\qquad {\bf x} \in B_i^{(n)} \\
    n\, \dist({\bf x}, \partial B_i ) \qquad\,\, {\bf x}\in B_i \setminus B_i^{(n)} \\
    0 \qquad\qquad\qquad\qquad {\bf x}\notin B_i
    \end{cases} . $$
    Put $ f \coloneqq \sum_{i=1}^N f_i^{(n)}$, which is continuous since $ f_i^{(n)} $'s are continuous by construction. Moreover, we compute
    $$ \int_{B_i} f = \int_{B_i} f_i^{(n)} \geq |B_i^{(n)}|.$$
    Hence
    \begin{align*}
        \Vert \Pi \Vert_{\mathrm{op}} \geq \Vert \Pi f \Vert_{B} & = \sup_{B \in \mathcal{B}^n (\Omega) } \frac{1}{|B|} \left\vert \int_B \sum_{i=1}^N \left( \int_{B_i} f \right) \ell_{B_i} \right\vert \\
        &= \sup_{B \in \mathcal{B}^n (\Omega) } \frac{1}{|B|} \left\vert \sum_{i=1}^N \left( \int_{B_i} f \right) \int_B \ell_{B_i} \right\vert  \\  
        &\geq \sup_{B \in \mathcal{B}^n (\Omega) } \frac{1}{|B|} \left\vert \sum_{i=1}^N |B_i^{(n)} | \int_B \ell_{B_i} \right\vert \\
        &= \sup_{B \in \mathcal{B}^n (\Omega) } \frac{1}{|B|} \left\vert \sum_{i=1}^N \frac{|B_i|}{{|B_i|}}|B_i^{(n)} | \int_B \ell_{B_i} \right\vert .
    \end{align*}
    Now, observe that
    $$ \frac{|B_i^{(n)}|}{|B_i|} = \left( 1 - \frac{1}{r_i n} \right) \geq \left( 1 - \frac{1}{r' n} \right) ,$$
    being $ r' \coloneqq \max_i r_i $. Thus 
     \begin{align*}
        & \sup_{B \in \mathcal{B}^n (\Omega) } \frac{1}{|B|} \left\vert \sum_{i=1}^N \frac{|B_i|}{{|B_i|}}|B_i^{(n)} | \int_B \ell_{B_i} \right\vert 
         \geq \left( 1 - \frac{1}{r' n} \right) \sup_{B \in \mathcal{B}^n (\Omega) } \frac{1}{|B|} \left\vert \sum_{i=1}^N |B_i| \int_B \ell_{B_i} \right\vert ,
    \end{align*}
    whence the claim, letting $ n \to \infty $.
\end{proof}

The above proof immediately extends to weaker regularity on $ f $. Moreover, if the hypothesis of disjointness $ B_i \cap B_j = \emptyset $ is not satisfied, the quantity \eqref{eq:newLeb} is nevertheless an upperbound for the norm of the interpolation operator.

\begin{proposition} \label{prop:operatornormbound}
    Let $ \{ B_i \}_{i=1}^N $ (possibly not disjoint) be unisolvent for $ \P_d(\R^n) $.
    Then
    $$ \Vert \Pi \Vert_{\mathrm{op}} \leq \Lambda_d .$$
\end{proposition}

\begin{proof}
    Untangling definitions, we have
    \begin{align*}
        \Vert \Pi \Vert_{\mathrm{op}} & = \sup_{\Vert f \Vert_B = 1} \Vert \Pi f \Vert_B = \sup_{\Vert f \Vert_B = 1} \left\Vert \sum_{i=1}^N \left( \int_{B_i} f \right) \ell_{B_i} \right\Vert_B \\
        & \leq \sup_{\Vert f \Vert_B = 1} \sup_{B \in \mathcal{B}^n} \frac{1}{|B|} \sum_{i=1}^N \left| \int_B \left( \int_{B_i} f \right) \ell_{B_i} \right| \\
        & = \sup_{\Vert f \Vert_B = 1} \sup_{B \in \mathcal{B}^n} \frac{1}{|B|} \sum_{i=1}^N \frac{|B_i|}{|B_i|} \left| \int_B \left( \int_{B_i} f \right) \ell_{B_i} \right| \\
        & = \sup_{\Vert f \Vert_B = 1} \sup_{B \in \mathcal{B}^n} \frac{1}{|B|} \sum_{i=1}^N \frac{\left|\int_{B_i} f \right|}{|B_i|} |B_i | \left| \int_B  \ell_{B_i} \right| \\
        & \leq \sup_{\Vert f \Vert_B = 1} \sup_{B \in \mathcal{B}^n} \frac{1}{|B|} \sum_{i=1}^N \Vert f \Vert_B |B_i | \left| \int_B  \ell_{B_i} \right| \\
        & = \sup_{B \in \mathcal{B}^n} \frac{1}{|B|} \sum_{i=1}^N |B_i | \left| \int_B  \ell_{B_i} \right| = \Lambda_d .
    \end{align*}
    The claim is proved.
\end{proof}

Since $ \Pi $ is also a projector, as proved in Proposition \ref{prop:projector}, the following bound can be established.

\begin{proposition} \label{prop:extensionclassicalestimate}
    Let $ \{ B_i \}_{i=1}^N $ be unisolvent for $ \P_d(\R^n) $, and let $ p^* $ be such that $ \Vert f - p^* \Vert_B \leq \Vert f - p \Vert_B $ for each $ p \in \P_d(\R^n) $. Then
    $$ \Vert f - \Pi f \Vert_B \leq (1+\Lambda) \Vert f - p^* \Vert_B .$$
\end{proposition}

\begin{proof}
    We compute 
    \begin{align*}
        \Vert f - \Pi f \Vert_B & = \Vert f - p^* + p^* - \Pi f \Vert_B \leq \Vert f - p^* \Vert_B + \Vert \Pi f - p^* \Vert_B \\
        & = \Vert f - p^* \Vert_B + \Vert \Pi f - \Pi p^* \Vert_B = \Vert f - p^* \Vert_B + \Vert \Pi (f - p^*) \Vert_B \\
        & \leq \Vert f - p^* \Vert_B + \Vert \Pi \Vert_{\mathrm{op}} \Vert f - p^* \Vert_B = \left( 1 + \Vert \Pi \Vert_{\mathrm{op}} \right) \Vert f - p^* \Vert_B.
    \end{align*}
    The claim now follows from Proposition \ref{prop:operatornormbound}.
\end{proof}

A crucial aspect for a meaningful definition of a Lebesgue constant \eqref{eq:newLeb} is that it must be affected only by the mutual positioning of the supports $ B_i $ (with respect to the reference domain). The following collection of results show that \eqref{eq:newLeb} satisfies this request. We first offer a simplification.

\begin{lemma} \label{lem:simplifiedLeb}
One has
\begin{equation} \label{eq:SimplifiedLeb}
    \Lambda_d = \sup_{B \in \mathcal{B}^n (\Omega)} \frac{1}{|R|^n} \sum_{i=1}^N |r_i|^n \left| \int_B \ell_{B_i} \right| ,
\end{equation}
with $ R $ being the radius of $ B $ and $ r_i $ that of $ B_i $.
\end{lemma}

\begin{proof}
    It suffices to plug \eqref{eq:ballvolume} into \eqref{eq:newLeb} to obtain the simplification, as
    \begin{align*}  
    \Lambda_d & = \sup_{B \in \mathcal{B}^n (\Omega)} \frac{1}{|B|} \sum_{i=1}^N | B_i | \left\vert \int_B \ell_{B_i} \right\vert = \sup_{B \in \mathcal{B}^n (\Omega)} \frac{\Gamma(\frac{n}{2}+1)}{\pi^n R^n} \sum_{i=1}^N \frac{\pi^n r_i^n}{\Gamma(\frac{n}{2}+1)} \left\vert \int_B \ell_{B_i} \right\vert \\ & = \sup_{B \in \mathcal{B}^n (\Omega)} \frac{1}{|R|^n} \sum_{i=1}^N |r_i|^n \left| \int_B \ell_{B_i} \right|,
    \end{align*}
    so \eqref{eq:SimplifiedLeb} is established.
\end{proof}

We show that \eqref{eq:newLeb} is invariant under translations, multiples of the identity and the action of $ SO(n) $; the Lebesgue constant is thus invariant under mappings that preserve $n$-balls. First of all, notice that a combination of the above kinds of mappings gives an affinity of $ \R^n $
\begin{equation} \label{eq:varphi}
    \varphi: \quad {\bf x} \mapsto \lambda A {\bf x} + {\bf b},
\end{equation}
with $ \lambda \in \R $, $ b \in \R^n $ and $ A \in SO(n) $. If $ {\bf p}$, $ {\bf q }$ are points in $ \R^n $, one has
$$ \varphi ({\bf p}) - \varphi ({\bf q}) = \lambda A {\bf p} + b - \lambda A {\bf q} - b = \lambda A ({\bf p} - {\bf q}) ,$$
whence $ \Vert \varphi ({\bf p}) - \varphi ({\bf q}) \Vert_2 = | \lambda | \Vert {\bf p}-{\bf q} \Vert_2 $ since $ A \in SO(n) $.
This adjusts volume terms in \eqref{eq:SimplifiedLeb}, as taking $ {\bf p } $ on the boundary of $ B $ and $ {\bf q } $ as the centre of $ B $, one has $ r = \Vert {\bf p}-{\bf q} \Vert_2 $, whence
\begin{equation} \label{eq:transfradius}
\Vert \varphi({\bf p})- \varphi({\bf q}) \Vert_2 = |\lambda| r .
\end{equation}

\begin{proposition} \label{prop:invarianceLeb}
    Let $ \varphi $ be as in \eqref{eq:varphi}
    and let
    $$ \Lambda_d^{\varphi} \coloneqq \sup_{\varphi(B) \in \mathcal{B}^n (\Omega)}  \frac{1}{|\varphi(B)|} \sum_{i=1}^N | \varphi(B_i) | \left \vert \int_{\varphi(B)} \ell_{\varphi(B_i)} \right \vert $$
    be the Lebesgue constant \eqref{eq:newLeb} associated with $ \left\{ \varphi (B_i) \right\}_{i=1}^N $. Then
    $$ \Lambda_d = \Lambda_d^{\varphi} .$$
\end{proposition}

\begin{proof}
First, we shall understand how Lagrange bases $ \left\{\ell_{\varphi(B_i)} \right\}_{i=1}^N$ for $ \left\{ \varphi (B_i) \right\}_{i=1}^N $ are done. By the change of variable, we have
$$ \delta_{i,j} = \int_{B_j} \ell_{B_i} = \int_{\varphi(B_j)} \ell_{B_i} \circ \varphi^{-1},$$
and we put $ \ell_{\varphi(B_i)} \coloneqq \ell_{B_i} \circ \varphi^{-1} $. Since $ \varphi^{-1} $ is a non degenerate affinity, $ \left\{ \ell_{\varphi(B_i)} \right\}_{i=1}^N $ is a basis for polynomials of degree $ d $ which satisfies the duality constraint \eqref{eq:dualityLB}. As a consequence, it is the Lagrange basis for $ \left\{\varphi(B_i) \right\}_{i=1}^N$. As a consequence, we also get that
    $$ \int_{B} \ell_{B_i} = \int_{\varphi(B)} \ell_{\varphi(B_i)}$$
    for each $ B \in \mathcal{B}^n (\Omega) $ and $ i = 1, \ldots, N $.
    To conclude the proof, it is now sufficient apply Lemma \ref{lem:simplifiedLeb} to the definition of $ \Lambda_d^\varphi $ and plug Eq. \eqref{eq:transfradius} in to obtain
    \begin{align*}
        \Lambda_d^\varphi & = \sup_{\varphi(B) \in \Omega} \frac{1}{||\lambda| R|^n} \sum_{i=1}^N ||\lambda| r_i|^n \left| \int_{\varphi(B)} \ell_{\varphi(B_i)} \right| \\
        & = \sup_{\varphi(B) \in \Omega} \frac{1}{|R|^n} \sum_{i=1}^N |r_i|^n \left| \int_{\varphi(B)} \ell_{\varphi(B_i)} \right| \\
        & = \sup_{B \in \Omega} \frac{1}{|R|^n} \sum_{i=1}^N |r_i|^n \left| \int_B \ell_{B_i} \right| = \Lambda_d,
    \end{align*}
    and the proof is concluded.
\end{proof}

Proposition \ref{prop:invarianceLeb} guarantees that Eq. \eqref{eq:newLeb} defines a generalised Lebesgue constant in the sense of \cite{ARR20}. As a consequence, \eqref{eq:newLeb} proves to be a meaningful generalisation of the nodal case to this setting.


\section{A specialisation to $ \mathbb{R}^2 $} \label{sect:r2}

In this section we make concrete the results of Section \ref{sect:theory} in the two-dimensional setting. 
The identification of unisolvent points is supported by a large bibliography, from the classical \cite{CY77, GM82} to the more recent \cite{BSV12, GSV11}. Theorem \ref{thm:generalunisolvence} puts us in a position to bridge several approaches adopted in literature: indeed, it does not distinguish between unisolvent points that are computed on a simplex, on a polygon or on a sphere, thus substantially enriching the possible choices of centres.

The case of orbits has a neat outline in this case.  Let $ B^{2} $ be the disc of radius $ 1 $ centered at $ 0 $. First of all, let us show that Theorem \ref{thm:generalunisolvence} and Theorem \ref{thm:unisolvenceorbits} coincide when $ d \leq 1 $ (the only non trivial situation is $ d = 1 $). Indeed $ \dim \P_1 (\R^2) = 3 $ and any three distinct nodes $ \{ \xi_1, \xi_2, \xi_3 \} $ are unisolvent for $\P_1 (\R^2)$. Theorem \ref{thm:generalunisolvence} grants that for any $ r > 0 $ the discs $ \{ B_{\xi_1,r}, B_{\xi_2,r}, B_{\xi_3,r} \} $. On the other hand, these three points define a unique circumference, and such an orbit is saturated with respect to $\P_1 (\R^2)$. Hence Theorem \ref{thm:unisolvenceorbits} yields the same conclusion. 

We evidence the separation between even and odd cases. As example, let us fix $ d = 2 $. Since $ \dim \mathbb{P}_2 (\R^2) = 6 $, we aim at placing $ 6 $ discs with prescribed radius $ \rho $. Fix any distance $ 0< r_1 < r - \rho $ and consider the orbit $ S_{r_1}^{1} $. The number of centres we shall place on $ S_{r_1}^{1} $ is prescribed by Eq. \eqref{eq:dimrestr} which we recall that, for $n=2$, it corresponds to $2d+1$. By substituting $ n = 2 $ and $ d = 2 $, we get
$$ \dim \mathbb{P}_2 ( S_{r_1}^{1} ) = 5 .$$
In fact, once we fix any $ B_{\xi_1,r} = B ({\xi_1 \cdot e^{i \theta},r} ) $, by passing to polar coordinates, it is easy to see that the mapping
$$
\theta \mapsto q(\theta) \coloneqq \int_{B ({\xi_1 \cdot e^{i \theta},r} )}p(x,y) \de x \de y
$$
defines a trigonometric polynomial of degree $ 2 $ in one variable $ \theta $. Since any distinct $ 5 $ points $ 0 \leq \theta_1 < \ldots < \theta_5 < 2 \pi $ are unisolvent for such a space, it is sufficient to take distincts centres. 
Moreover, by applying Lemma \ref{lem:factorisationSr}, one gets that once one place, not in the same orbit, any ball centered at $ {\bf c} \in B^2 \setminus S_r^{1} $, we get the unisolvence for the whole space $ \P_2 ( \R^2) $. Indeed, it follows by Theorem \ref{thm:generalunisolvence} that balls of fixed radius centered at such points are unisolvent as well.

Let now $ d = 3 $, hence $ \dim \mathbb{P}_3 (\R^2) = 10 $. Recalling Eq. \ref{eq:dimrestr} we obtain that $ \dim (\mathbb{P}_3 (S_{r_1}^1)) = 7 $. As a consequence, we may place seven points on an orbit specified by $ r_1 $ and the unisolvence implies that either $ p(x,y) = 0 $ or $ p(x,y) $ factorises through the equation $ x^2+y^2=r_1^2$.
We are left with the placement of the resulting three points. \\
Saturating another orbit specified by $ 0 < r_2 < r-\rho$, $r_2 \neq r_1 $ thus grants global unisolvence, see the above case for $ d = 1 $. By iterating this reasoning one easily reconstructs the proof of Theorem \ref{thm:unisolvenceorbits} when $ n = 2 $.

\begin{remark}
    In accordance with the nodal case analysed in \cite{Bos91}, it is worth pointing out that in the even cases the degenerate orbit given by only one point is required. That this point can be \emph{any} point in $ B $ (not lying on any already considered orbits) is a consequence of Theorem \ref{thm:unisolvenceorbits} and shows the necessity of requiring saturated orbits. 
\end{remark}
\begin{remark}
    When $ n = 2 $, term $ d^* $ appearing in Theorem \ref{thm:unisolvenceorbits} is easily computed to be $ d^* = \lceil \frac{d}{2} + 1 \rceil $.
\end{remark}

\section{Some remarkable fails} \label{sect:fails}

In Section \ref{sect:theory} we have exhibited several ways of obtaining unisolvent collections of balls. By a density argument, one may further deduce that almost every randomly picked set $ \{ B_i \}_{i=1}^N $, with $ N = \dim \P_d (\R^n) $, is unisolvent for $ \dim \P_d (\R^n) $. It is thus of some interest to study if there are some peculiar sets that are \emph{not} unisolvent. We dedicate this section to a couple of interesting families that fail this scope, working examples out in the two-dimensional setting; however, we remark that they can be easily extended to greater dimension. 

\begin{figure}[!h]
\centering
\begin{tikzpicture}
  \draw[->] (-3, 0) -- (3.5, 0) node[right] {$x$};
  \draw[->] (0, -3) -- (0, 3.5) node[above] {$y$};

\fill [color={rgb,255:red,79; green,120; blue,153}, opacity=0.6,even odd rule] (0,0) circle[radius=0.4cm];
\fill [color={rgb,255:red,15; green,186; blue,189}, opacity=0.6,even odd rule] (0,0) circle[radius=0.4cm] circle[radius=0.9cm];
\fill [color={rgb,255:red,33; green,103; blue,94}, opacity=0.6,even odd rule] (0,0) circle[radius=0.9cm] circle[radius=1.2cm];
\fill [color={rgb,255:red,251; green,189; blue,60}, opacity=0.6,even odd rule] (0,0) circle[radius=1.2cm] circle[radius=1.5cm];
\fill [color={rgb,255:red,188; green,19; blue,31}, opacity=0.6,even odd rule] (0,0) circle[radius=1.5cm] circle[radius=2.2cm];
\fill [color={rgb,255:red,140; green,16; blue,68}, opacity=0.6,even odd rule] (0,0) circle[radius=2.2cm] circle[radius=2.5cm];
  
\draw[thin,dashed] (0,0) circle (0.4cm);
\draw[thin,dashed] (0,0) circle (0.9cm);
\draw[thin,dashed] (0,0) circle (1.2cm);
\draw[thin,dashed] (0,0) circle (1.5cm);
\draw[thin,dashed] (0,0) circle (2.2cm);

\draw[thin] (0,0) circle (2.5 cm);
\end{tikzpicture} \quad
\begin{tikzpicture}
  \draw[->] (-3, 0) -- (3.5, 0) node[right] {$x$};
  \draw[->] (0, -3) -- (0, 3.5) node[above] {$y$};

\filldraw[color={rgb,255:red,33; green,103; blue,94}, opacity=0.6] (0.3,0) circle (0.2cm);
\filldraw[color={rgb,255:red,33; green,103; blue,94}, opacity=0.6] (-0.3,0) circle (0.2cm);
\filldraw[thin,dashed,color={rgb,255:red,251; green,189; blue,60}, opacity=0.6] (1,0) circle (0.2cm);
\filldraw[thin,dashed,color={rgb,255:red,251; green,189; blue,60}, opacity=0.6] (-1,0) circle (0.2cm);
\filldraw[thin,dashed,color={rgb,255:red,188; green,19; blue,31}, opacity=0.6] (2,0) circle (0.2cm);
\filldraw[thin,dashed,color={rgb,255:red,188; green,19; blue,31}, opacity=0.6] (-2,0) circle (0.2cm);

\draw[thin] (0,0) circle (2.5 cm);

\draw[thin,dotted, color={rgb,255:red,33; green,103; blue,94}, opacity=1] (0,0) circle (0.3 cm);
\draw[thin,dotted,color={rgb,255:red,251; green,189; blue,60}, opacity=1] (0,0) circle (1 cm);
\draw[thin,dotted,color={rgb,255:red,188; green,19; blue,31}, opacity=1] (0,0) circle (2 cm);

\end{tikzpicture}
\caption{Two non-unisolvent configurations. Left: a non-unisolvent configuration with prescribed centre. Right: a non-unisolvent configuration with not saturated orbits.} \label{fig:annulus}
\end{figure}
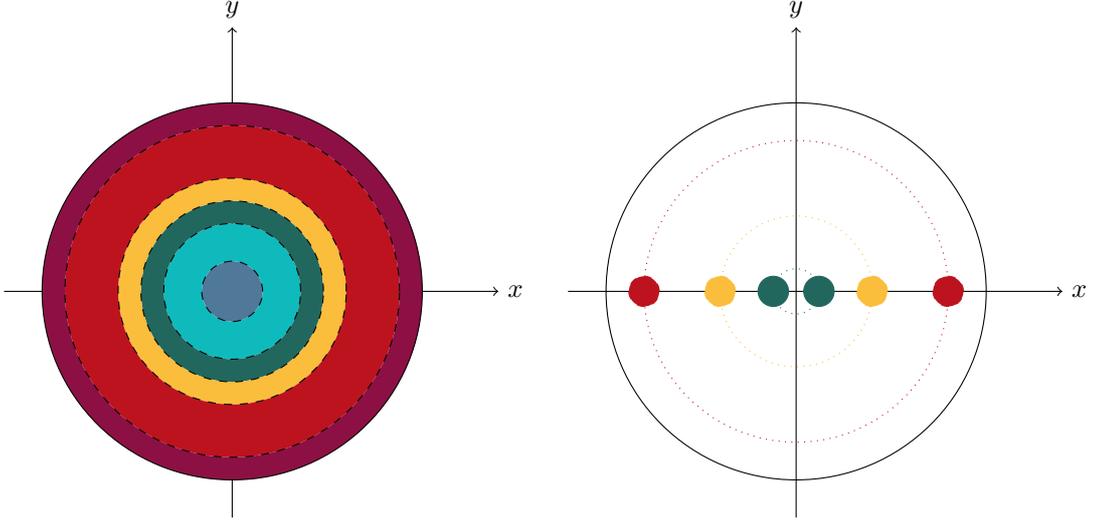

\subsection{Fixed centre, changing radius}

This example shows that it is better not to exceed with symmetry.

Suppose to have data on annulus as in Figure \ref{fig:annulus}, left. By the linearity of the integral, this is equivalent to have same information on concentric discs. To check that this configuration of discs cannot be unisolvent for any $ d \geq 1 $, it suffices to show that block of the Vandermonde matrix
$$ V_{i,j} \coloneqq \int_{B_j} x^{i-1} y^{d-i+1} \de x \de y $$
corresponding to degree $ d $ homogeneous polynomials always presents an entire row of zeros. For any $ r_j > 0 $, let $ p(x,y) \in \P_d(\R^2) $. Exploiting polar coordinates, we compute
\begin{align*}
    \int_{B (0, r_j)} p(x,y) \de x \de y & = \int_0^{r_j} \int_0^{2\pi} r_j^{d+1} p(\cos(\theta), \sin(\theta) \de r \de \theta \\
    & = \frac{r_j^{d+2}}{d+2} \int_{0}^{2\pi} p(\cos(\theta),\sin(\theta)) \de \theta .
\end{align*}
It is immediate to observe that both the second and the third row are null. Indeed, in both cases we get
$$ V_{1,j} = \frac{r_j^{d+2}}{d+2} \int_{0}^{2\pi} \cos(\theta) \de \theta = 0 = \frac{r_j^{d+2}}{d+2} \int_{0}^{2\pi} \sin (\theta) \de \theta = V_{2,j} $$
for each $ r_j \in \R $, hence for any $ j = 1, \ldots, N $.

\subsection{Not saturated orbits} \label{sect:fail2}

Theorem \ref{thm:unisolvenceorbits} prescribed a sufficient condition for unisolvence based on orbits. In particular, following \cite[Lemma $3.1$]{Bos91}, we obtained an upper bound of the centres that we may place on any orbit. We show the convenience of saturating orbits.

We construct a case of a large number of non saturated orbits where at least two rows of the Vandermonde matrix are proportional. This is depicted in Figure \ref{fig:annulus} on the right. Let us fix a disc $ B_1 $, say the leftmost one, we, then, observe that any other disc $ B_j $ is obtained as $ B_j = B_1 + (\xi_j, 0) $. Therefore we have that
$$ \int_{B_j} p(x,y) \de x \de y = \int_{B_1} p(x-\xi_j,y) \de x \de y .$$
If we consider again the monomial basis, we have that the row
$$ \int_{B_j} y^t \de x \de y= \int_{B_1} y^t \de x \de y = c_t $$
associated with $ p(x,y) = y^t $ is constant for each $ j = 1, \ldots, N $. As a consequence, they are all multiple of the first row, which is
$$ V_{1,j} = \int_{B_j} 1 \de x \de y = |B_j| = \pi r_j^2. $$
It follows that the Vandermonte matrix is degenerate and hence the collection $ \{ B_j\}_{j=1}^N $ is not unisolvent for $ \P_d (\R^2) $.

\section{Numerical experiments} \label{sect:numeric}

Among others, the aim of this section is to provide a comparison between the point-set approach of Section \ref{sect:unisolvencepoints} and the orbit approach of Section \ref{sect:unisolvenceorbits}. We will show that they are both valuable and that some fails may occur also on the numerical level. In Section \ref{sect:interpolation} we produce interpolation tests. In Section \ref{sect:numLeb} some Lebesgue constants are estimated in this framework. Matlab demos and codes are available at
\url{https://github.com/gelefant/InterpDISC}.

\subsection{A convenient basis}

\begin{figure}[!h]
   \centering
        {\includegraphics[height=2.20in]{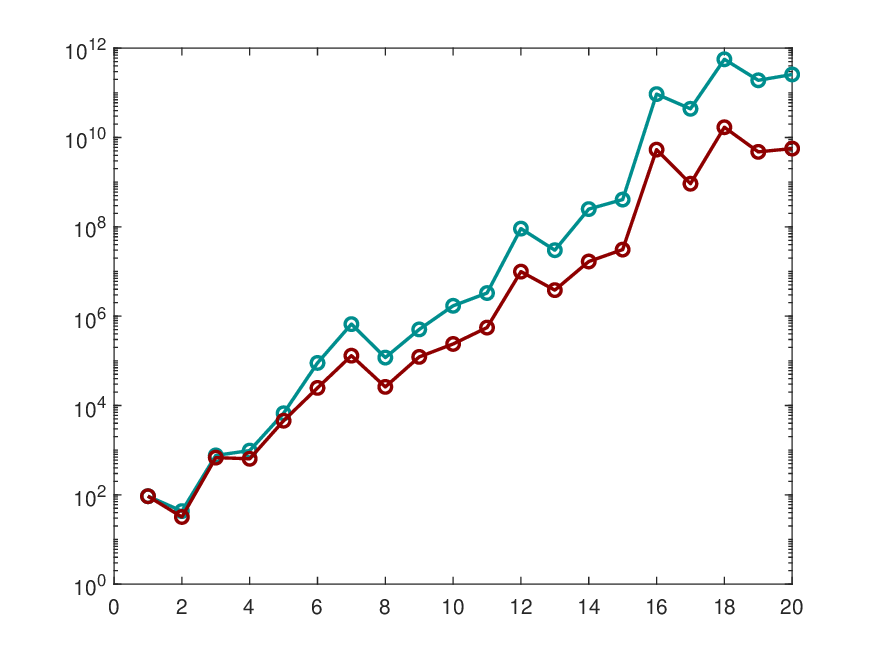}}
     \quad
        {\includegraphics[height=2.20in]{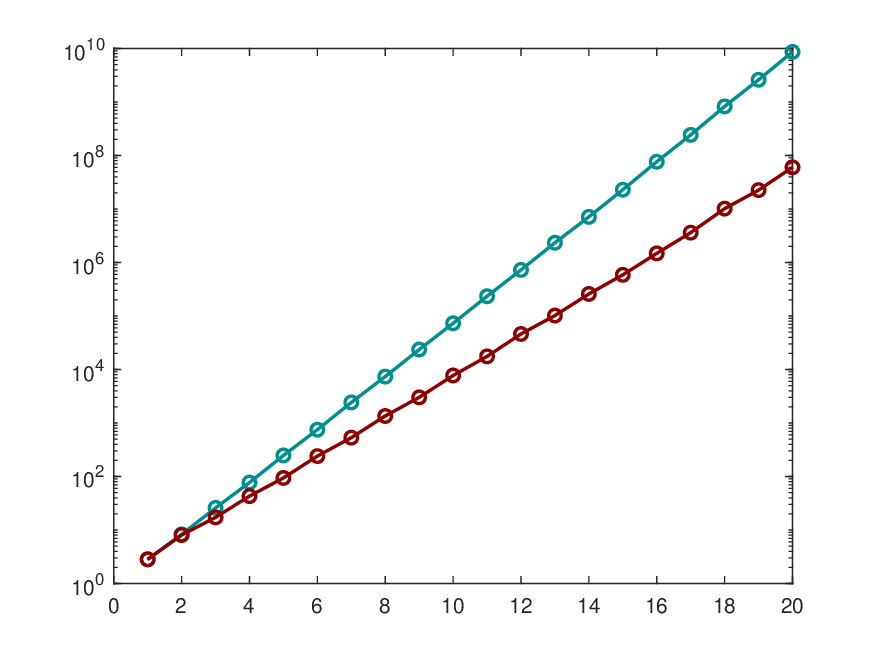}}     
\caption{Conditioning of the Vandermonde matrix: Chebyshev polynomials (red) vs monomials (blue). Halton points (left) and Chebyshev orbits (right) are depicted.} \label{fig:cond}
\end{figure}

For all our purposes, we shall need at hand explicit expressions for the Lagrange dual basis $ \ell_{B_i} $ defined in \eqref{eq:dualityLB}. Since there are not, in general, closed formulae, we may fix a convenient basis $ p_i $, compute the relative Vandermonde matrix and invert it to obtain the coefficients of \eqref{eq:LagrBases} (see, e.g., \cite[Chapter 3]{BruniThesis}). The choice of the convenient basis is thus driven by two aspects: the existence of exact quadrature formulae on the disc and the well conditioning of the resulting Vandermonde matrix. Two natural candidates are thus the monomial and the Chebyshev basis. We thus tested these two bases on different collections of discs, based on both point sets (also quasi-random) and orbits. As we can notice in Figure \ref{fig:cond}, we have that using Chebyshev basis leads in general to a better conditioned Vandermonde matrix. In the following we thus stick with this choice.

\begin{nat}
In the rest of the paper, we say that a collection of balls is based on some point sets if their centres belong to this set. We thus refer to:
\begin{itemize}
    \item Chebyshev orbits: orbits whose radius is a non-negative Chebyshev points in $[-1,1]$. Discs whose centres belong to this points are named \emph{Chebyshev discs}.
    \item equidistant orbits: orbits whose radius is of the form $ r_j = \frac{j}{N}$, for $ 0 \leq j < N $;
    \item Halton discs: discs whose centres belong to the Halton sequence in two dimension (see, e.g., \cite{Niederreiter92}). We also refer to these discs as \emph{quasi-random};
    \item optimal (and quasi-optimal) discs: discs centered at nodes that present low nodal Lebesgue constants on the disc. For these sets we refer to \cite{MS19}. We slightly push all the nodes towards the centre $ {\bf 0 }$ in order to not have points lying on the boundary and construct on such points discs with small radius.
    \end{itemize}
    In the case of orbits, the centres on each orbit are placed to span equidistant arcs. Rotations of orbits seem not to affect sensibly the problem, as already observed in the nodal case \cite{BSV12}. 
\end{nat}

\begin{figure}[!h]
   \centering
        {\includegraphics[height=2.20in]{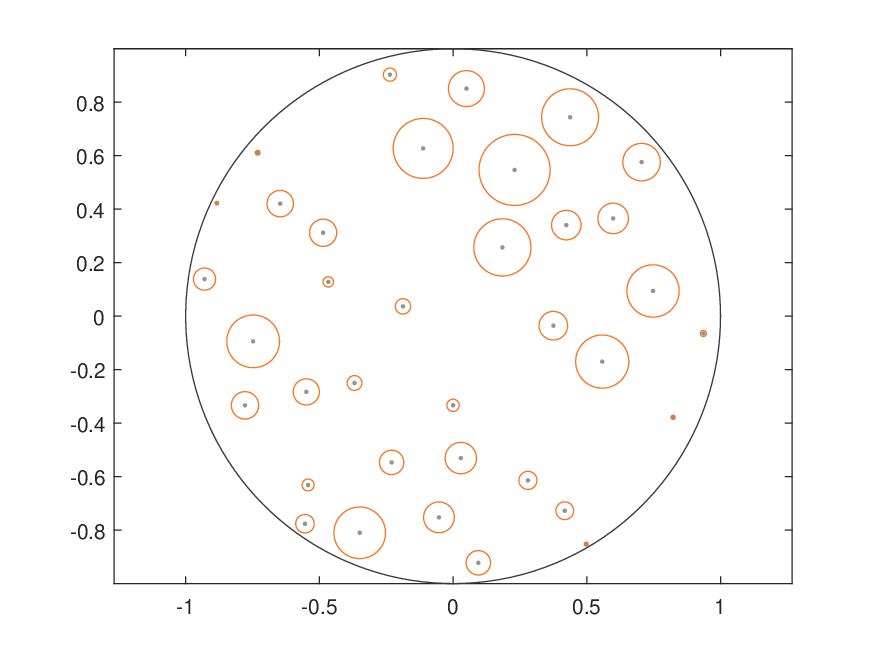}}
    \quad 
        {\includegraphics[height=2.20in]{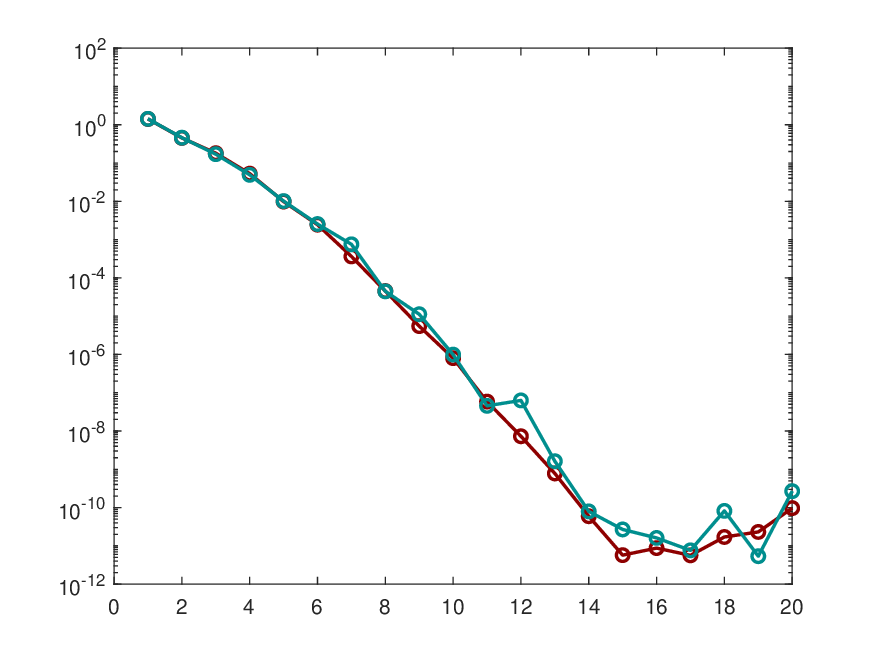}}     
\caption{Left: discs with random radii placed on Halton discs. Right: interpolation errors when random radii are taken into account on Chebyshev orbits, both with (blue) and without (red) overlapping.}
\label{fig:NoOrb}
\end{figure}

This choice is motivated by the numerics. In fact, as we shall see in the forthcoming pages, the radius of the supporting balls seem not particularly relevant in determining the numerical quality of the interpolation. For an account of this, we address to Figure \ref{fig:NoOrb} on the right. Such a figure depicts an interpolation test (detailed in the following section) where we constrained centres on Chebyshev orbits and let radii vary randomly also allowing overlapping. The corresponding interpolation errors (right panel) turn out to be not very sensitive to this choice. As a consequence, we will confine ourselves in the case where there is no intersection of the supports and the radii are fixed as the maximum possible, so that we are in fact estimating the norm of the interpolation operator.

\subsection{Interpolation} \label{sect:interpolation}

We test the efficiency of the interpolator \eqref{eq:definterp} by comparing different families of supports. In particular, we perform some interpolation tests with different functions and read the interpolation error as the supports of integration vary.

\begin{figure}[!h]
   \centering
        {\includegraphics[height=2.20in]{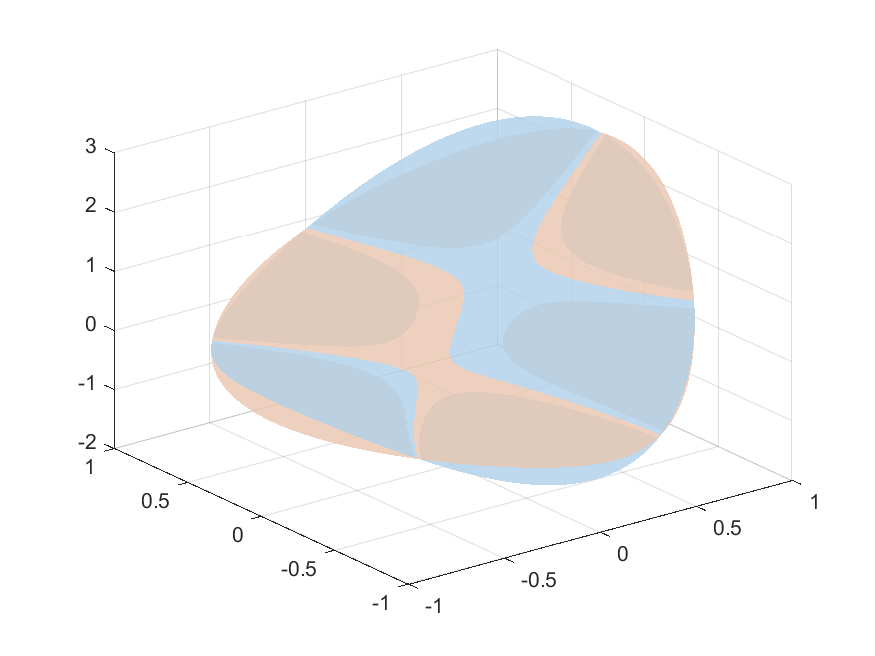}} 
\caption{The exact solution (blue) versus the interpolated (pink) for degree $ d = 2 $.}
\label{fig:interpf1}
\end{figure}

To begin with, we consider the function
$$ f_1 (x,y) = e^x \sin(x+y) .$$
In Figure \ref{fig:interpf1} we plot the graph of $ f_1 $ together to the graph of the interpolant of degree $ d = 2 $.

\begin{figure}[!h]
   \centering
        \includegraphics[height=2.20in]{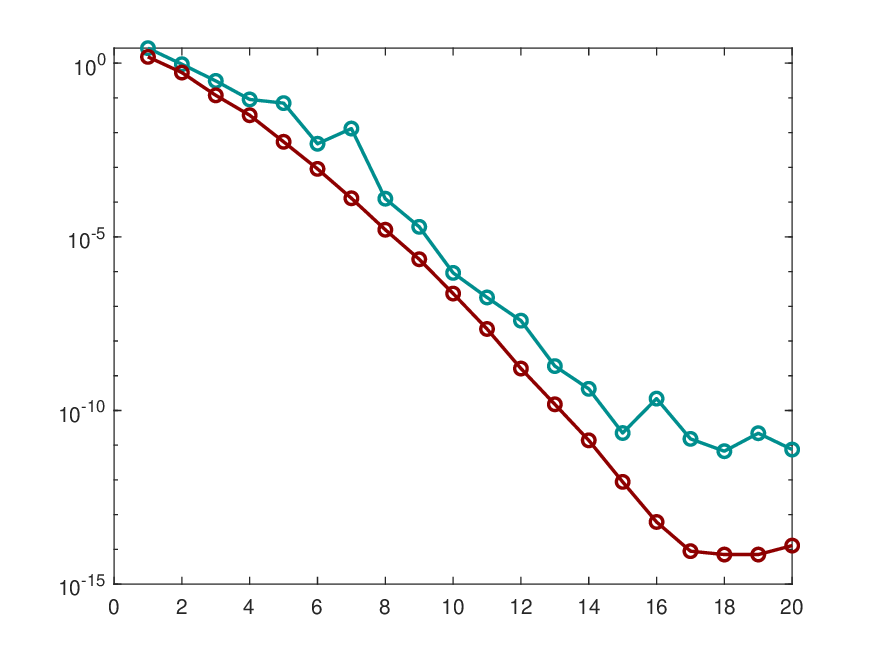}
        \includegraphics[height=2.20in]{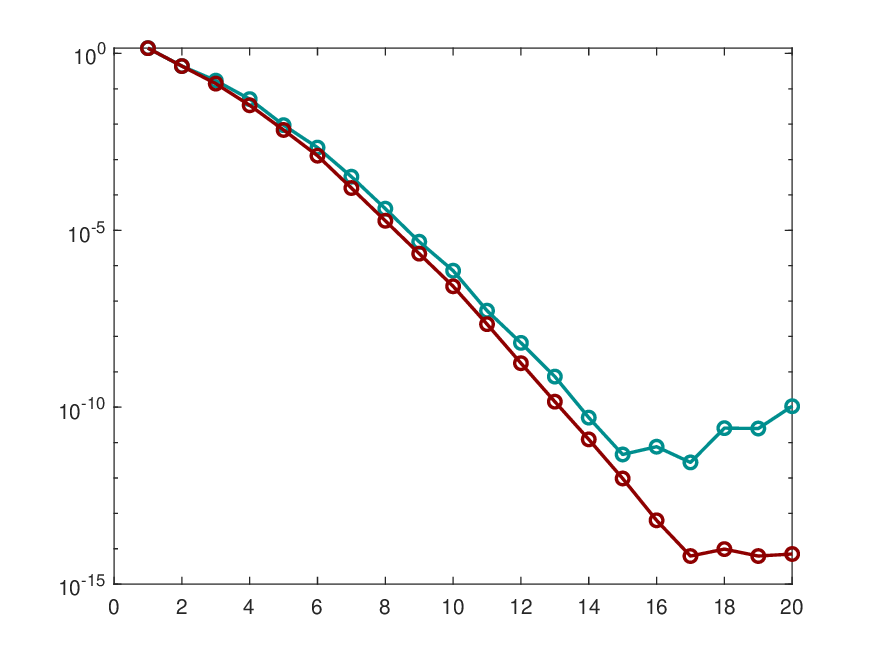}
\caption{Error $ \Vert f_1 - \Pi_d f_1 \Vert_0 $ versus degree $ d$. Left: Halton discs (blue) and optimal discs (red). Right: equidistant orbits (blue) and Chebyshev orbits (red). All the radii are constant and chosen to avoid overlappings.}
\label{fig:errf1}
\end{figure}

We let the total polynomial degree $ d $ increase and report $ \Vert f_1 - \Pi_d f_1 \Vert_0 $ in Figure \ref{fig:errf1}. As we see, both the point set case and the orbit case offer a solid decrease with the total polynomial degree $ d $.

\begin{figure}[!h]
   \centering
        {\includegraphics[height=2.20in]{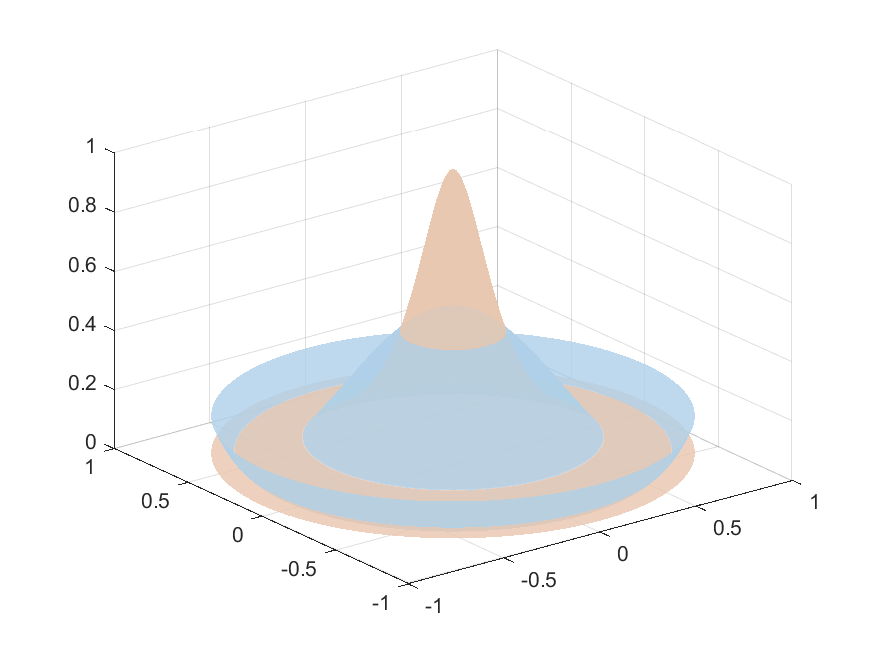}} \quad
        {\includegraphics[height=2.20in]{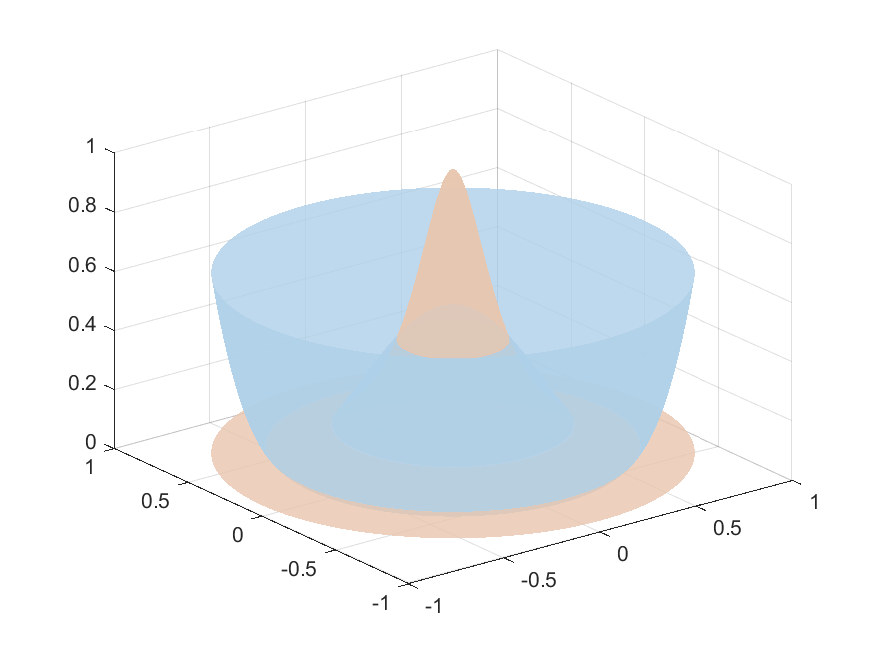}}     
\caption{The exact function $f_2$ and two interpolated. One based on discs with radius driven by a uniform distribution (right), the other with radius based on a Chebyshev distribution (left), for $ d = 5 $.}
\label{fig:interpRunge}
\end{figure}

As Runge-like phenomena associated with interpolators \eqref{eq:definterp} have been observed, both in the one dimensional \cite{BE23} and the multivariate \cite{BruniRunge} framework, we propose a function that is expected to show such a behaviour. It is based on the construction of \cite{RothThesis} and reads
$$ f_2 (x,y) = \frac{1}{25(x^2+y^2)+1} .$$
In this case the choice of supports seriously affects the reliability of the interpolated, as depicted in Figure \ref{fig:interpRunge}. In accordance with the expectations, the interpolated captures the function only away from the boundary, where instability becomes evident.

\begin{figure}[!h]
   \centering
        \includegraphics[height=2.20in]{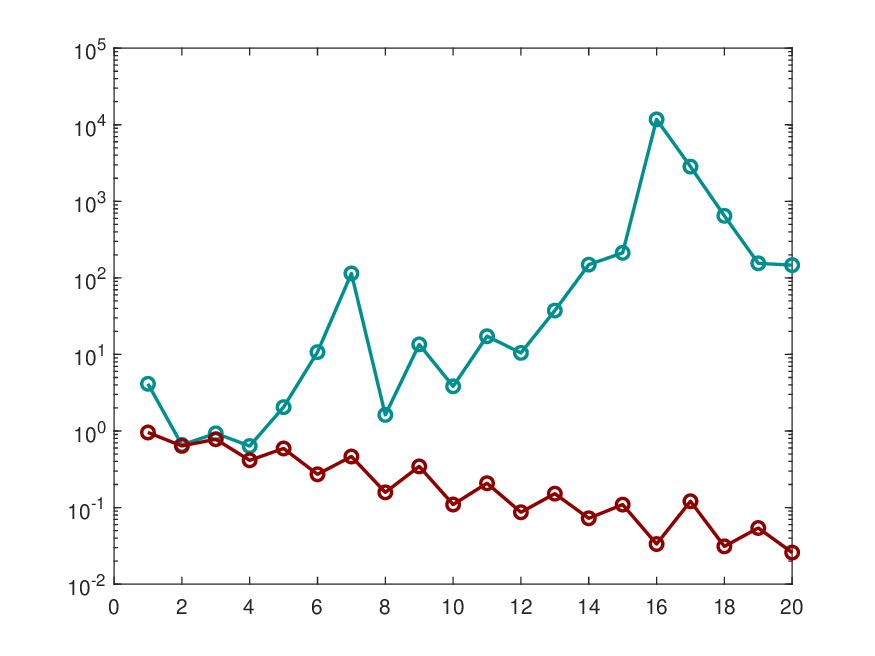}
        \includegraphics[height=2.20in]{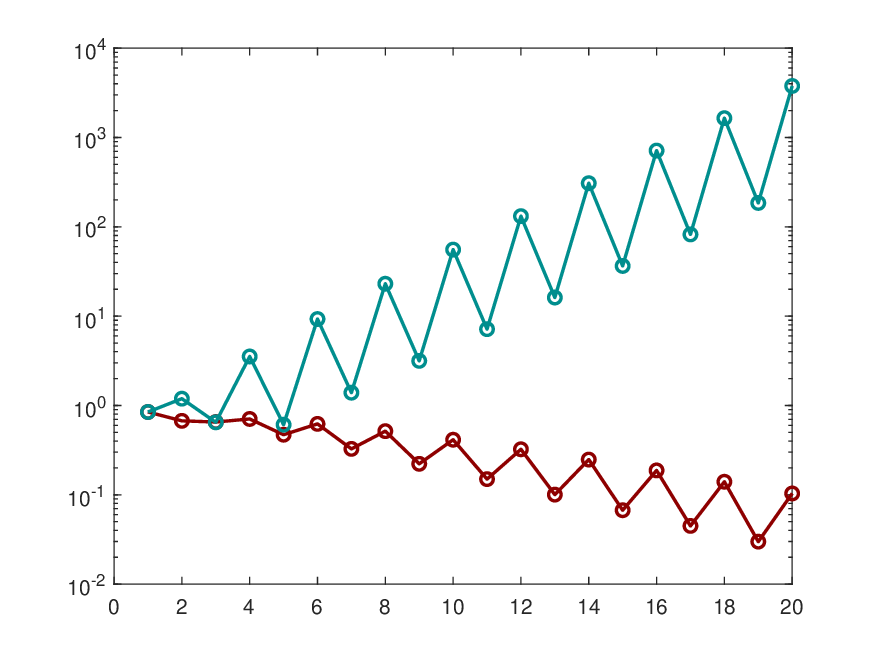}
\caption{Error $ \Vert f_2 - \Pi_d f_2 \Vert_0 $ versus degree $ d$. Left: Halton discs (blue) and optimal discs (red). Right: equidistant orbits (blue) and Chebyshev orbits (red). All the radii are constant and chosen to avoid overlappings.}
\label{fig:errf2}
\end{figure}

The trend of the errors $ \Vert f_2 - \Pi f_2 \Vert_0 $ with respect to the total degree $ d $ is reported in Figure \ref{fig:errf2}. In such a figure we show that we are able to find both meaningful and poor point sets (left hand panel) as well as significant and bad orbits (right hand panel). This confirms that increasing the total degree of the interpolated is not a solution, providing an extension of the Runge phenomenon to this framework. As a consequence, the choice of appropriate discs turn out to be a relevant task for this problem.

Motivated by this test and Proposition \ref{prop:extensionclassicalestimate}, we thus look for discs that offer a low Lebesgue constant in the sense of \eqref{eq:newLeb}. To this pursuit we dedicate the next section.

\subsection{Estimation of Lebesgue constants} \label{sect:numLeb}

In order to estimate the Lebesgue constant $ \Vert \Pi \Vert_{\mathrm{op}}$, we compute \eqref{eq:newLeb} over a rich collection of balls $ \{D_k\}_{k=1}^M $ of radius $ \rho_k $. We point out that, even if supports $ B_i $'s intersect, Proposition \ref{prop:operatornormbound} ensures that \eqref{eq:newLeb} overestimates the norm of the interpolation operator; as a conseuquence, the capability of containing this quantity is as relevant as having precise estimates for $ \Vert \Pi \Vert_{\mathrm{op}}$. On the other hand, since $ \Pi $ is a projector onto polynomials of degree $d$, its norm is bounded from below by
\begin{equation} \label{eq:lowerbound}
    \Vert \Pi \Vert_{\mathrm{op}} \geq c_n \cdot \begin{cases}
        \log d, \quad \text{if } n = 1, \\
        d^\frac{n-1}{2}, \quad \text{if } n > 1,
    \end{cases}
\end{equation}
as established in \cite{Sundermann}.

Let us recall and expand \eqref{eq:newLeb} as
\begin{align*}
    \Lambda_d & = \sup_{D_k} \frac{1}{|{D_k}|} \sum_{i=1}^N |B_i| \left\vert \int_{D_k} \ell_{B_i} \right\vert = \sup_{D_k} \frac{1}{|{D_k}|} \sum_{i=1}^N |B_i| \left\vert \int_{D_k} \sum_{j=1}^n V^{-1}_{i,j} p_i \right\vert \\ &  = \sup_{D_k} \frac{1}{|{D_k}|} \sum_{i=1}^N |B_i| \left\vert \sum_{j=1}^n V^{-1}_{i,j} \int_{D_k} p_i \right\vert 
     = \sup_{D_k} \frac{1}{|{\rho_k}|} \sum_{i=1}^N |r_i| \left\vert \sum_{j=1}^n V^{-1}_{i,j} \int_{D_k} p_i \right\vert,
\end{align*}
where $ \rho_k $ is the radius of $ D_k $. Hence, define the $ M \times N $ matrix
$$ W_{i,k} \coloneqq \int_{D_i} p_k $$
and the square matrices
$$ \boldsymbol{\rho} = \begin{pmatrix} |\rho_1|^{-d} & 0 & \ldots & 0 \\
0 & \ddots & \ddots & 0 \\
0 & \ddots & \ddots & 0 \\
0 & \ldots & 0 & |\rho_M|^{-d} 
\end{pmatrix} \quad \text{and} \quad
{\bf r} = \begin{pmatrix} |r_1|^d & 0 & \ldots & 0 \\
0 & \ddots & \ddots & 0 \\
0 & \ddots & \ddots & 0 \\
0 & \ldots & 0 & |r_N|^d
\end{pmatrix} .
$$
It follows that, for each $ k $,
$$ \frac{1}{|{D_k}|} \sum_{i=1}^N |B_i| \left\vert \int_{D_k} \ell_{B_i} \right\vert = \sum_{t=1}^N \left\vert \left(\boldsymbol{\rho} W V^{-T} {\bf r}\right)_{k,t} \right\vert ,$$
whence
$$ \Lambda_d = \Vert \boldsymbol{\rho} W V^{-T} {\bf r} \Vert_1, $$
being $ \Vert \cdot \Vert_1 $ the $1$-norm for matrices (i.e. the maximum by row-wise sum of absolute values). Again, we point out that all the quantities at play are exactly computed, provided that radii are.

\begin{figure}[!h]
   \centering
    \subfigure[]
        {\includegraphics[height=2.20in]{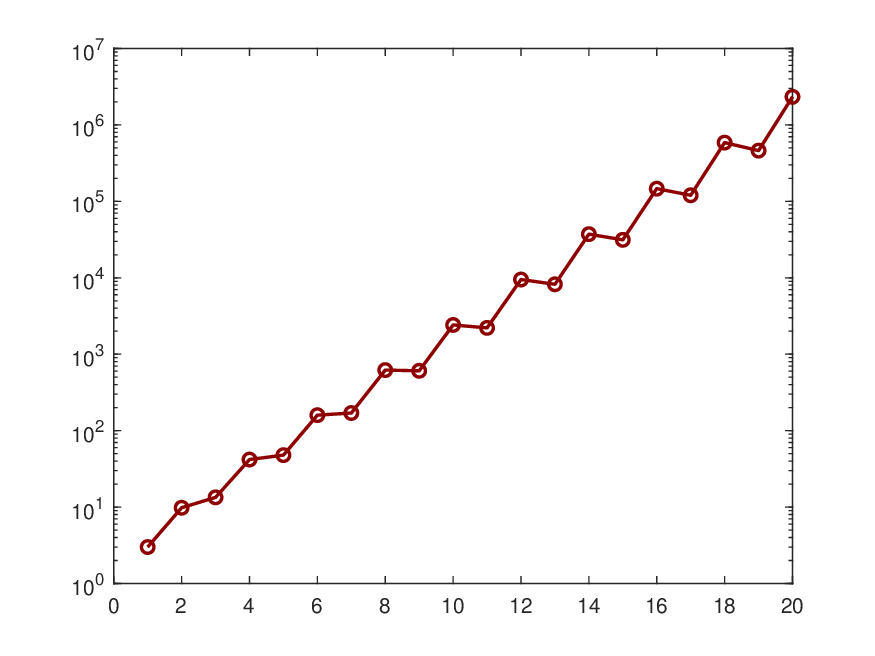}}
     ~ 
    \subfigure[]
        {\includegraphics[height=2.20in]{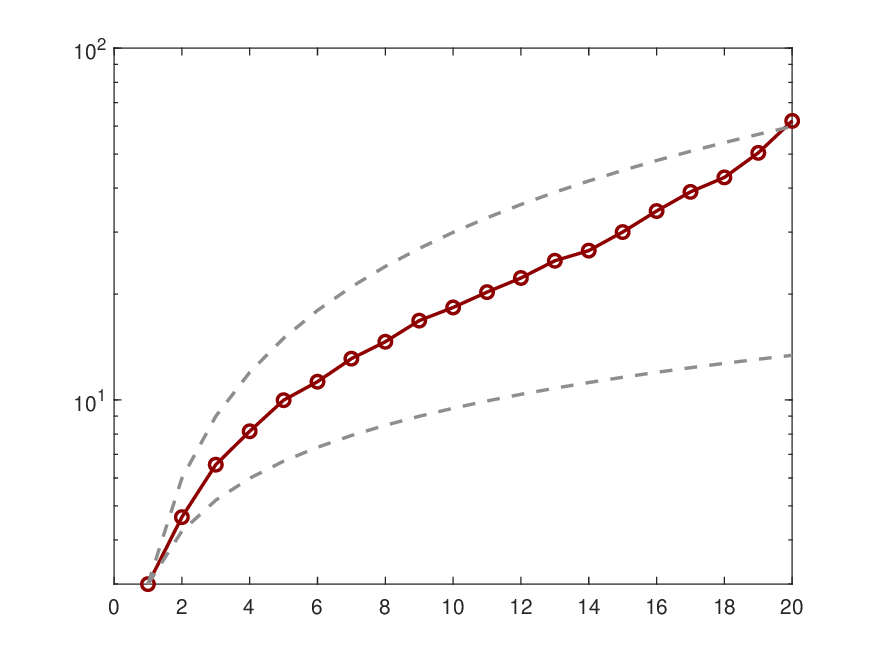}}     
             ~ 
    \subfigure[]
        {\includegraphics[height=2.20in]{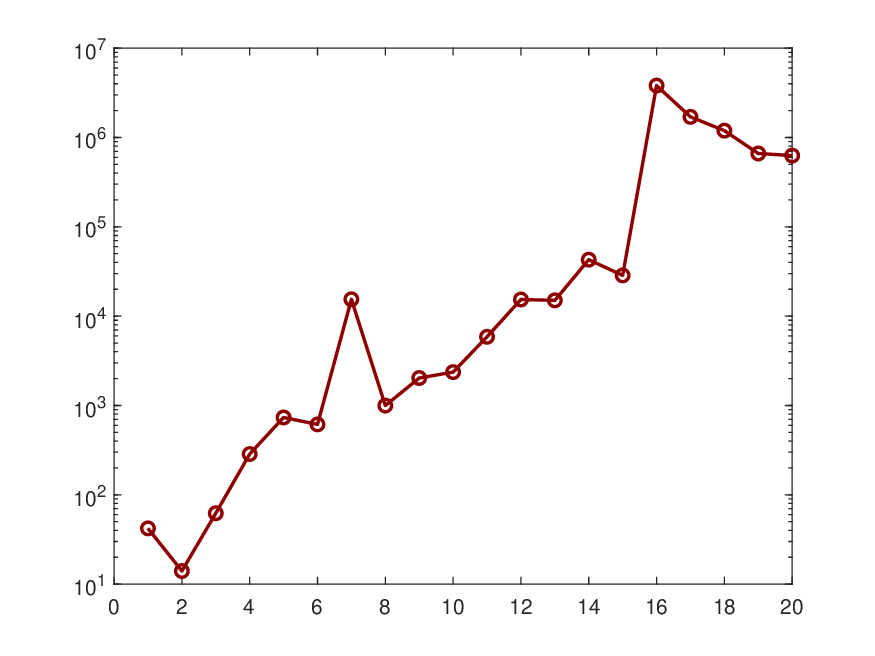}}  
             ~ 
    \subfigure[]
        {\includegraphics[height=2.20in]{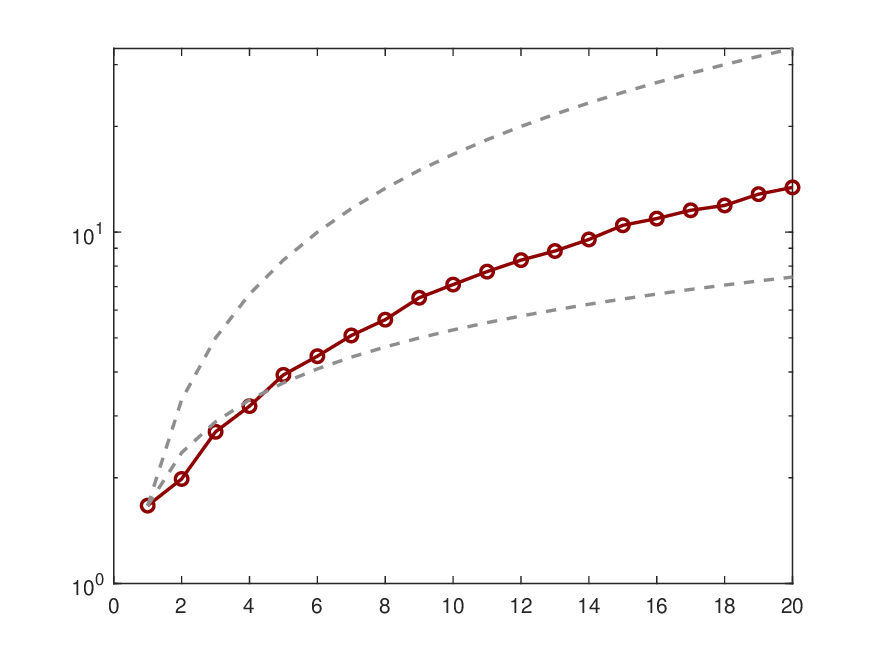}}  
\caption{Lebesgue constants for equidistributed orbits (a), Chebyshev orbits (b), Halton discs (c) and optimal discs (d) as the total degree $ d $ increases.}
\label{fig:leb}
\end{figure}

In Figure \ref{fig:leb} we compare several choices of supports in terms of their Lebesgue constant \eqref{eq:newLeb}. The panel (a) depicts the Lebesgue constant associated with equidistributed orbits. The trend is clearly exponential and this is consistent with the right hand panel of Figure \ref{fig:errf2}, which shows poor interpolation properties associated with this choice. Likewise does the panel (c), that depicts the behaviour of Halton discs. The left hand panel of Figure \ref{fig:errf2} relates this trend with a undesirable behaviour of the corresponding interpolator.

In panels (b) and (d), results are compared with a linear bound (the upper one, dashed grey) and a sublinear one ($ \sqrt{d} $, taking for ease $ c_N = 1 $, dashed grey), which is prescribed for optimality by \eqref{eq:lowerbound}. Chebyshev orbits seem to offer interesting performances, although they struggle a bit as $ d $ approaches $ 20 $; see panel (b). In contrast, discs that are centred at nodes that present a low (nodal) Lebesgue constant seem to be quite close to the optimal bound, as depicted in panel (d).

\section{Conclusions} \label{sect:conclusions}

In this paper we have proposed and analysed an interpolation strategy which is based on averages over regular domains such as $n$-balls. We provided sufficient conditions for unisolvence by transforming integrals into evaluations of a convenient auxiliary polynomial. Two strategy then came at play: the first relies on the interpolation of functions on points sets, whereas the second on orbits induced by the action of $ SO(n) $ on the unit $n$-ball. In both cases, we showed that the radius of the support balls is not relevant for unisolvence, although the case of orbits allows for more various choices. To compare them, we introduced a concept of Lebesgue constant which in fact measures or (upper) bounds the norm of the resulting interpolation operator. Numerical experiments showed that the choice of radii is in fact irrelevant to stability features. Several point sets- and orbit-based collections of points have been analysed numerically, evidencing advantages and drawbacks of each of them. Numerical results confirm the effectiveness of the theoretical side.

\section*{Acknowledgements}

This research has been accomplished within Rete ITaliana di Approssimazione (RITA), the thematic group on \lq\lq Teoria dell'Ap\-pros\-si\-ma\-zio\-ne e Applicazioni\rq\rq $\,$(TAA) of the Italian Mathematical Union and partially supported by GNCS-IN$\delta$AM. The first author is funded by IN$\delta$AM and supported by Universit\`a di Padova.

\printbibliography

@article {ChR16,
    AUTHOR = {Christiansen, S. H. and Rapetti, F.},
     TITLE = {On high order finite element spaces of differential forms},
   JOURNAL = {Math. Comp.},
  FJOURNAL = {Mathematics of Computation},
    VOLUME = {85},
      YEAR = {2016},
    NUMBER = {298},
     PAGES = {517--548},
}

@misc{BE23,
      title={Polynomial Interpolation of Function Averages on Interval Segments}, 
      author={Bruni Bruno, L. and Erb, W.},
      year={2023},
      eprint={2309.00328},
      archivePrefix={arXiv},
      primaryClass={math.NA}
}

@article {BBZ22,
    AUTHOR = {Bruni Bruno, L. and Zampa, E.},
     TITLE = {Unisolvent and minimal physical degrees of freedom for the
              second family of polynomial differential forms},
   JOURNAL = {ESAIM Math. Model. Numer. Anal.},
  FJOURNAL = {ESAIM. Mathematical Modelling and Numerical Analysis},
    VOLUME = {56},
      YEAR = {2022},
    NUMBER = {6},
     PAGES = {2239--2253},
}

@phdthesis{BruniThesis,
	author       = {Bruni Bruno, L.}, 
	title        = {Weights as degrees of freedoom for high order Whitney finite elements},
	school       = {University of Trento},
	year         = {2022},
	note = {Available at: \url{https://theses.hal.science/tel-04067201/}}
}

@article {Bos91,
    AUTHOR = {Bos, L.},
     TITLE = {On certain configurations of points in {${\bf R}^n$} which are
              unisolvent for polynomial interpolation},
   JOURNAL = {J. Approx. Theory},
  FJOURNAL = {Journal of Approximation Theory},
    VOLUME = {64},
      YEAR = {1991},
    NUMBER = {3},
     PAGES = {271--280},
}

@article {GSV11,
    AUTHOR = {Gentile, M. and Sommariva, A. and Vianello, M.},
     TITLE = {Polynomial interpolation and cubature over polygons},
   JOURNAL = {J. Comput. Appl. Math.},
  FJOURNAL = {Journal of Computational and Applied Mathematics},
    VOLUME = {235},
      YEAR = {2011},
    NUMBER = {17},
     PAGES = {5232--5239},
}

@article {BSV12,
    AUTHOR = {Briani, M. and Sommariva, A. and Vianello, M.},
     TITLE = {Computing {F}ekete and {L}ebesgue points: simplex, square,
              disk},
   JOURNAL = {J. Comput. Appl. Math.},
  FJOURNAL = {Journal of Computational and Applied Mathematics},
    VOLUME = {236},
      YEAR = {2012},
    NUMBER = {9},
     PAGES = {2477--2486},
}

@article {ARR20,
    AUTHOR = {Alonso Rodr\'{\i}guez, A. and Rapetti, F.},
     TITLE = {On a generalization of the {L}ebesgue's constant},
   JOURNAL = {J. Comput. Phys.},
  FJOURNAL = {Journal of Computational Physics},
    VOLUME = {428},
      YEAR = {2021},
     PAGES = {Paper No. 109964, 4},
}

@article {CY77,
    AUTHOR = {Chung, K. C. and Yao, T. H.},
     TITLE = {On lattices admitting unique {L}agrange interpolations},
   JOURNAL = {SIAM J. Numer. Anal.},
  FJOURNAL = {SIAM Journal on Numerical Analysis},
    VOLUME = {14},
      YEAR = {1977},
    NUMBER = {4},
     PAGES = {735--743},
}

@article {GM82,
    AUTHOR = {Gasca, M. and Maeztu, J. I.},
     TITLE = {On {L}agrange and {H}ermite interpolation in {${\bf
              R}\sp{k}$}},
   JOURNAL = {Numer. Math.},
  FJOURNAL = {Numerische Mathematik},
    VOLUME = {39},
      YEAR = {1982},
    NUMBER = {1},
     PAGES = {1--14},
}

@article {Phung21,
    AUTHOR = {Phung, V. M.},
     TITLE = {Hermite interpolation on the unit sphere and limits of
              {L}agrange projectors},
   JOURNAL = {IMA J. Numer. Anal.},
  FJOURNAL = {IMA Journal of Numerical Analysis},
    VOLUME = {41},
      YEAR = {2021},
    NUMBER = {2},
     PAGES = {1441--1464},
}

@article {Phung17,
    AUTHOR = {Phung, V. M.},
     TITLE = {Polynomial interpolation in {$\mathbb{R}^2$} and on the unit
              sphere in {$\mathbb{R}^3$}},
   JOURNAL = {Acta Math. Hungar.},
  FJOURNAL = {Acta Mathematica Hungarica},
    VOLUME = {153},
      YEAR = {2017},
    NUMBER = {2},
     PAGES = {289--317},
}

@article {Hesthaven98,
    AUTHOR = {Hesthaven, J. S.},
     TITLE = {From electrostatics to almost optimal nodal sets for
              polynomial interpolation in a simplex},
   JOURNAL = {SIAM J. Numer. Anal.},
  FJOURNAL = {SIAM Journal on Numerical Analysis},
    VOLUME = {35},
      YEAR = {1998},
    NUMBER = {2},
     PAGES = {655--676},
}

@incollection {HarrisonLeb,
    AUTHOR = {Harrison, J.},
     TITLE = {Continuity of the integral as a function of the domain},
      NOTE = {Dedicated to the memory of Fred Almgren},
   JOURNAL = {J. Geom. Anal.},
  FJOURNAL = {The Journal of Geometric Analysis},
    VOLUME = {8},
      YEAR = {1998},
    NUMBER = {5},
     PAGES = {769--795},
}

@article {ABR20,
    AUTHOR = {Alonso Rodr\'{\i}guez, A. and Bruni Bruno, L. and Rapetti,
              F.},
     TITLE = {Towards nonuniform distributions of unisolvent weights for
              high-order {W}hitney edge elements},
   JOURNAL = {Calcolo},
  FJOURNAL = {Calcolo. A Quarterly on Numerical Analysis and Theory of
              Computation},
    VOLUME = {59},
      YEAR = {2022},
    NUMBER = {4},
     PAGES = {Paper No. 37, 29},
}

@article {Rapetti07,
    AUTHOR = {Rapetti, F.},
     TITLE = {High order edge elements on simplicial meshes},
   JOURNAL = {M2AN Math. Model. Numer. Anal.},
  FJOURNAL = {M2AN. Mathematical Modelling and Numerical Analysis},
    VOLUME = {41},
      YEAR = {2007},
    NUMBER = {6},
     PAGES = {1001--1020},
}

@incollection {Xu04,
    AUTHOR = {Xu, Y.},
     TITLE = {Polynomial interpolation on the unit sphere and on the unit ball},
      NOTE = {Approximation and applications},
   JOURNAL = {Adv. Comput. Math.},
  FJOURNAL = {Advances in Computational Mathematics},
    VOLUME = {20},
      YEAR = {2004},
    NUMBER = {1-3},
     PAGES = {247--260},
}

@article {Xu03,
    AUTHOR = {Xu, Y.},
     TITLE = {Polynomial interpolation on the unit sphere},
   JOURNAL = {SIAM J. Numer. Anal.},
  FJOURNAL = {SIAM Journal on Numerical Analysis},
    VOLUME = {41},
      YEAR = {2003},
    NUMBER = {2},
     PAGES = {751--766},
}

@incollection {Bojanov06,
    AUTHOR = {Bojanov, B.},
     TITLE = {Interpolation and integration based on averaged values},
 BOOKTITLE = {Approximation and probability},
    SERIES = {Banach Center Publ.},
    VOLUME = {72},
     PAGES = {25--47},
 PUBLISHER = {Polish Acad. Sci. Inst. Math., Warsaw},
      YEAR = {2006},
}

@article {Warburton06,
    AUTHOR = {Warburton, T.},
     TITLE = {An explicit construction of interpolation nodes on the
              simplex},
   JOURNAL = {J. Engrg. Math.},
  FJOURNAL = {Journal of Engineering Mathematics},
    VOLUME = {56},
      YEAR = {2006},
    NUMBER = {3},
     PAGES = {247--262},
}

@article {MS19,
    AUTHOR = {Meurant, G. and Sommariva, A.},
     TITLE = {On the computation of sets of points with low {L}ebesgue
              constant on the unit disk},
   JOURNAL = {J. Comput. Appl. Math.},
  FJOURNAL = {Journal of Computational and Applied Mathematics},
    VOLUME = {345},
      YEAR = {2019},
     PAGES = {388--404},
}

@article {BE21,
    AUTHOR = {Berrut, J.-P. and Elefante, G.},
     TITLE = {Bounding the {L}ebesgue constant for a barycentric rational
              trigonometric interpolant at periodic well-spaced nodes},
   JOURNAL = {J. Comput. Appl. Math.},
  FJOURNAL = {Journal of Computational and Applied Mathematics},
    VOLUME = {398},
      YEAR = {2021},
     PAGES = {113664, 11},
}

@article {TFR22,
    AUTHOR = {Takaki, N. and Forbes, G. W. and Rolland, J. P.},
     TITLE = {Schemes for cubature over the unit disk found via numerical
              optimization},
   JOURNAL = {J. Comput. Appl. Math.},
  FJOURNAL = {Journal of Computational and Applied Mathematics},
    VOLUME = {407},
      YEAR = {2022},
     PAGES = {114076, 19},
}

@article {CK00,
    AUTHOR = {Cools, R. and Kim, K. J.},
     TITLE = {A survey of known and new cubature formulas for the unit disk},
   JOURNAL = {Korean J. Comput. Appl. Math.},
  FJOURNAL = {The Korean Journal of Computational \& Applied Mathematics. An
              International Journal},
    VOLUME = {7},
      YEAR = {2000},
    NUMBER = {3},
     PAGES = {477--485},
}

@article {BruniRunge,
    AUTHOR = {Alonso Rodr\'{\i}guez, A. and Bruni Bruno, L. and Rapetti,
              F.},
     TITLE = {Whitney edge elements and the {R}unge phenomenon},
   JOURNAL = {J. Comput. Appl. Math.},
  FJOURNAL = {Journal of Computational and Applied Mathematics},
    VOLUME = {427},
      YEAR = {2023},
     PAGES = {115117, 9},
}

@phdthesis{RothThesis,
	title    = {Nodal configurations and Voronoi tessellations for triangular spectral elements},
	school   = {University of Victoria},
	author   = {Roth, M. J.},
	year     = {2005},
}

@article {BuhmannIntegration,
    AUTHOR = {Buhmann, M. D. and Dai, F. and Niu, Y.},
     TITLE = {Discretization of integrals on compact metric measure spaces},
   JOURNAL = {Adv. Math.},
  FJOURNAL = {Advances in Mathematics},
    VOLUME = {381},
      YEAR = {2021},
     PAGES = {107602, 32},
}

@book{buhmann_2003, 
    place={Cambridge}, 
    series={Cambridge Monographs on Applied and Computational Mathematics}, 
    title={Radial Basis Functions: Theory and Implementations}, 
    publisher={Cambridge University Press}, 
    author={Buhmann, M. D.}, 
    year={2003},}

@article {Sundermann,
    AUTHOR = {S\"{u}ndermann, B.},
     TITLE = {On projection constants of polynomial spaces on the unit ball
              in several variables},
   JOURNAL = {Math. Z.},
  FJOURNAL = {Mathematische Zeitschrift},
    VOLUME = {188},
      YEAR = {1984},
    NUMBER = {1},
     PAGES = {111--117},
}

@book{Niederreiter92,
author = {Niederreiter, H.},
title = {Random Number Generation and Quasi-Monte Carlo Methods},
publisher = {Society for Industrial and Applied Mathematics},
year = {1992},
}

@article {Ibrahimoglu16,
    AUTHOR = {Ibrahimoglu, B. A.},
     TITLE = {Lebesgue functions and {L}ebesgue constants in polynomial
              interpolation},
   JOURNAL = {J. Inequal. Appl.},
  FJOURNAL = {Journal of Inequalities and Applications},
      YEAR = {2016},
     PAGES = {93, 15},
}

@article {DEM21,
    AUTHOR = {De Marchi, S. and Elefante, G. and Marchetti,
              F.},
     TITLE = {On {$(\beta,\gamma)$}-{C}hebyshev functions and points of the
              interval},
   JOURNAL = {J. Approx. Theory},
  FJOURNAL = {Journal of Approximation Theory},
    VOLUME = {271},
      YEAR = {2021},
     PAGES = {105634, 17},
}

@incollection {ABR23,
    AUTHOR = {Alonso Rodr\'{\i}guez, Ana and Bruni Bruno, Ludovico and Rapetti,
              Francesca},
     TITLE = {Flexible weights for high order face based finite element
              interpolation},
 BOOKTITLE = {Spectral and high order methods for partial differential
              equations {ICOSAHOM} 2020+1},
    SERIES = {Lect. Notes Comput. Sci. Eng.},
    VOLUME = {137},
     PAGES = {117--128},
 PUBLISHER = {Springer, Cham},
      YEAR = {2023},
}

\end{document}